\documentclass[graybox]{svmult}
\PassOptionsToPackage{bookmarksdepth=4}{hyperref}

	\usepackage{multicol}            
	\usepackage[bottom]{footmisc}    
	\usepackage{newtxtext}           
	\usepackage[mathscr]{euscript}
	\usepackage[explicit]{titlesec}  
	\let\I\relax                     
	\usepackage[%
		noload={caption,subcaption}, 
		eqreset=section,             
		thmreset=section,            
	]{myPreamble}
	\hypersetup{hypertexnames=false}
\AtBeginDocument{\renewcommand\thesection{\arabic{section}}}

	\makeatletter
		\providecommand*{\toclevel@titlech}{10}
		\edef\toclevel@authorch{\the\numexpr\toclevel@titlech}
	\makeatother

\makeatletter
	\def\chap@hangfrom#1{}%
\makeatother

	\AtBeginEnvironment{appendix}{%
		\renewcommand\appendixname{Appendix}%
		\setcounter{section}{0}%
		\renewcommand\thesection{\Alph{section}}%
		\titleformat{\section}{\normalfont\secsize\secstyle}{}{0em}{%
			Appendix \thesection\ifstrempty{#1}{}{. #1}%
		}%
	}

	\makeatletter
		\renewenvironment{svgraybox}{%
			\fboxsep=12pt\relax
			\begin{shaded}%
				\list{}{\leftmargin=\z@\topsep=\z@\rightmargin=\leftmargin\relax}%
				\expandafter\item\parindent=\svparindent
				\hskip-\listparindent
		}%
		{\endlist\end{shaded}}%
	\makeatother

	\makeatletter
		\renewcommand\alglinenumber[1]{\footnotesize\textbf{\thealgorithm}.\oldstylenums{\arabic{ALG@line}}:}
		\renewcommand\theALG@line{\thealgorithm.\oldstylenums{\arabic{ALG@line}}}
		\providecommand\theHALG@line{\thealgorithm.\oldstylenums{\arabic{ALG@line}}}
	\makeatother

	\newlist{proxexamples}{enumerate}{1}
	\AtBeginEnvironment{proxexamples}{%
		\let\olditem\item
		\renewcommand\item[1][]{%
			\def\itemtitle{#1}%
			\olditem
			\ifstrempty{#1}{}{%
				\phantomsection
				\addcontentsline{toc}{subsubsection}{\itemtitle}%
				\textbf{#1.}~~\ignorespaces%
			}%
		}
	}
	\setlist[proxexamples,1]{%
		label={\bf\alph*.},
		ref=\thesubsection\alph*,
		align=parleft,
	}
	\crefformat{proxexamplesi}{#2\S#1#3}
	\Crefformat{proxexamplesi}{#2\S#1#3}
	\crefname{proxexamplesi}{\S}{\S}
	\Crefname{proxexamplesi}{\S}{\S}

	\newcommand\refFBTN[1][ (\cref{alg:FBTN})]{%
		\hyperref[alg:FBTN]{{\sf FBTN}}#1%
	}%
	\newcommand\refCG[1][ (\cref{alg:CG0})]{%
		\hyperref[alg:CG0]{{\sf CG}}#1%
	}%

	\let\proj\relax
	\DeclareMathOperator\proj{\mathscr P}
	\DeclareMathOperator\sym{S}
	\DeclareMathOperator\eigs{eigs}
	\newcommand\FBE{\varphi_\gamma}
	\newcommandx\T[2][1=\gamma,2={}]{T_{#1}^{#2}}   
	\newcommandx\Res[2][1=\gamma,2={}]{R_{#1}^{#2}} 
	\newcommandx\lin[4][{1=\gamma,4={}}]{           
		\ifstrempty{#3}{%
			\ifstrempty{#2}{%
				\ell_{#1}^{#4}%
			}{%
				\ell_{#1}^{#4}\left(#2\right)
			}%
		}{%
			\ell_{#1}^{#4}\left(#2;\,#3\right)%
		}%
	}%

	\def\xMin{-1}
	\def\xMax{1}
	\def\yMin{-1}
	\def\yMax{1}
	\AtBeginEnvironment{axis}{%
		\pgfmathsetmacro\YMAXREAL{\yMax-(\yMax-\yMin)*0.05}
		\pgfmathsetmacro\YINFTY{\yMax-(\yMax-\yMin)*0.01}
	}
	\tikzset{%
		nosep/.style = {inner sep=0pt, outer sep=0pt},
		dash/.style  = {black, dotted},
		label/.style = {font=\scriptsize, inner sep=0.4pt, outer sep=1pt, preaction={fill=white, fill opacity=0.7}},
		dot/.style   = {mark=*, mark size=1.25pt, mark options={fill=white}},
	}
	\pgfplotsset{%
		ext_func/.style = {line width = 1pt},
		func/.style     = {ext_func, restrict y to domain={\yMin:\YMAXREAL}, samples=101},
		phi/.style      = {ext_func, color=phi},
		ME/.style       = {func, color=ME, dashed},
		FBE/.style      = {func, color=FBE},
		Q/.style        = {func, color=Q, line width = 0.3pt},
		myaxis/.style   = {%
			no markers,
			execute at begin axis = {\coordinate (origin) at (axis cs:0,0);},
			domain = \xMin:\xMax,
			xmin = \xMin,
			xmax = \xMax,
			ymin = \yMin,
			ymax = \yMax,
			xtick = \empty,
			ytick = \empty,
			axis lines = middle,
			enlargelimits = false,
			samples = 101,
			legend style = {%
				fill opacity = 0.95,
				font = \color{black}\scriptsize,
				legend cell align = left,
				rounded corners,
				anchor = north east,
				at = {(1,1)},
				line width = 0.4pt,
				outer sep = 0pt,
				inner sep = 1pt,
				cells = {line width = 1pt},
			},
		},
	}
	\colorlet{Q}{black!50}
	\colorlet{phi}{blue!80!black}
	\colorlet{ME}{black}
	\colorlet{FBE}{red!60!black}
	\newcommand\drawDot[2][]{
		\draw[dot, #1] plot coordinates {(#2)};
	}
	\newcommand\drawCoordX[2][]{%
		\draw[dash,#1] (origin -| #2) -- (#2);
	}
	\newcommand\drawCoordY[2][]{%
		\draw[dash,#1] (origin |- #2) -- (#2);
	}
	\newcommandx\drawCoord[3][{1={}},{2={}}]{%
		\drawCoordX[#1]{#3}
		\drawCoordY[#2]{#3}
	}
	\newcommand\xLabel[3][]{%
		\node[label, #1, anchor=north] at (#2 |- origin) {#3};
	}
	\newcommand\yLabel[3][]{%
		\node[label, #1, anchor=east] at (#2 -| origin) {#3};
	}
	\newcommand\plotExtFunc[3][]{
		\addplot[#1, func] {#2(x)};
		\addplot[#1, domain=#3, forget plot] {#2(x)};
	}

	\title{On the acceleration of forward-backward splitting via an inexact Newton method}
	\author{Andreas Themelis, Masoud Ahookhosh and Panagiotis Patrinos}
	\authorrunning{A. Themelis, M. Ahookhosh and P. Patrinos}
	\institute{%
		Andreas Themelis, Masoud Ahookhosh and Panagiotis Patrinos\at\TheAddressKU.
		This work was supported by the \emph{Research Foundation Flanders (FWO)} research projects G086518N and G086318N;
		\emph{KU Leuven internal funding} StG/15/043;
		\emph{Fonds de la Recherche Scientifique --- FNRS and the Fonds Wetenschappelijk Onderzoek --- Vlaanderen} under EOS Project no 30468160 (SeLMA).
		\\
		{\tt
			\{%
				\href{mailto:andreas.themelis@esat.kuleuven.be}{andreas.themelis},%
				\href{mailto:masoud.ahookhosh@esat.kuleuven.be}{masoud.ahookhosh},%
				\href{mailto:panos.patrinos@esat.kuleuven.be}{panos.patrinos}%
			\}%
			\href{mailto:andreas.themelis@esat.kuleuven.be,masoud.ahookhosh@esat.kuleuven.be,panos.patrinos@esat.kuleuven.be}{@esat.kuleuven.be}%
		}%
	}

\begin{document}

	\maketitle


	\abstract{%
		We propose a Forward-Backward Truncated-Newton method (FBTN) for minimizing the sum of two convex functions, one of which smooth.
		Unlike other proximal Newton methods, our approach does not involve the employment of variable metrics, but is rather based on a reformulation of the original problem as the unconstrained minimization of a continuously differentiable function, the \emph{forward-backward envelope (FBE)}.
		We introduce a generalized Hessian for the FBE that \emph{symmetrizes} the generalized Jacobian of the nonlinear system of equations representing the optimality conditions for the problem.
		This enables the employment of conjugate gradient method (CG) for efficiently solving the resulting (regularized) linear systems, which can be done inexactly.
		The employment of CG prevents the computation of full (generalized) Jacobians, as it requires only (generalized) directional derivatives.
		The resulting algorithm is globally (subsequentially) convergent, \(Q\)-linearly under an error bound condition, and up to \(Q\)-superlinearly and \(Q\)-quadratically under regularity assumptions at the possibly non-isolated limit point.%
	}

	


	\section{Introduction}
		In this work we focus on convex composite optimization problems of the form
\begin{equation}\label{eq:GenProb}
	\minimize_{x\in\R^n}\varphi(x)\equiv f(x)+g(x),
\end{equation}
where \(\func f{\R^n}{\R}\) is convex, twice continuously differentiable and with \(L_f\)-Lipschitz-continuous gradient, and \(\func{g}{\R^n}{\R\cup\set\infty}\) has a cheaply computable proximal mapping \cite{moreau1965proximite}.
To ease the notation, throughout the paper we indicate
\[
	\varphi_\star\coloneqq\inf\varphi
\quad\text{and}\quad
	\mathcal X_\star\coloneqq\argmin\varphi.
\]
Problems of the form \eqref{eq:GenProb} are abundant in many scientific areas such as control, signal processing, system identification, machine learning and image analysis, to name a few.
For example, when \(g\) is the indicator of a convex set then \eqref{eq:GenProb} becomes a constrained optimization problem, while for \(f(x)=\frac12\|Ax-b\|^2\) and \(g(x)=\lambda\|x\|_1\) it becomes the \(\ell_1\)-regularized least-squares problem (lasso) which is the main building block of compressed sensing.
When \(g\) is equal to the nuclear norm, then \eqref{eq:GenProb} models low-rank matrix recovery problems.
Finally, conic optimization problems such as linear, second-order cone, and semidefinite programs can be brought into the form of \eqref{eq:GenProb}, see \cite{lan2011primal}.

Perhaps the most well-known algorithm for problems in the form \eqref{eq:GenProb} is the forward-backward splitting (FBS) or proximal gradient method \cite{lions1979splitting,combettes2011proximal}, that interleaves gradient descent steps on the smooth function and \emph{proximal} steps on the nonsmooth one, see \cref{sec:Prox}.
Accelerated versions of FBS, based on the work of Nesterov \cite{nesterov2013gradient,beck2009fast,tseng2008accelerated}, have also gained popularity.
Although these algorithms share favorable global convergence rate estimates of order \(O(\varepsilon^{-1})\) or \(O(\varepsilon^{-1/2})\) (where \(\varepsilon\) is the solution accuracy), they are first-order methods and therefore usually effective at computing solutions of low or medium accuracy only.
An evident remedy is to include second-order information by replacing the Euclidean norm in the proximal mapping with that induced by the Hessian of \(f\) at \(x\) or some approximation of it, mimicking Newton or quasi-Newton methods for unconstrained problems \cite{becker2012quasi,lee2014proximal,lu2017randomized}.
However, a severe limitation of the approach is that, unless \(Q\) has a special structure, the computation of the proximal mapping becomes very hard.
For example, if \(\varphi\) models a lasso problem, the corresponding subproblem is as hard as the original problem.

In this work we follow a different approach by reformulating the nonsmooth constrained problem \eqref{eq:GenProb} into the smooth unconstrained minimization of the \emph{forward-backward envelope (FBE)} \cite{patrinos2013proximal}, a real-valued, continuously differentiable, exact penalty function for \(\varphi\).
Although the FBE might fail to be twice continuously differentiable, by using tools from nonsmooth analysis we show that one can design Newton-like methods to address its minimization, that achieve \(Q\)-superlinear asymptotic rates of convergence under nondegeneracy and (generalized) smoothness conditions on the set of solutions.
Furthermore, by suitably interleaving FBS and Newton-like iterations the proposed algorithm also enjoys good complexity guarantees provided by a global (non-asymptotic) convergence rate.
Unlike the approaches of \cite{becker2012quasi,lee2014proximal}, where the corresponding subproblems are expensive to solve, our algorithm only requires the inexact solution of a linear system to compute the Newton-type direction, which can be done efficiently with a memory-free CG method.

Our approach combines and extends ideas stemming from the literature on merit functions for \emph{variational inequalities} (VIs) and \emph{complementarity problems} (CPs), specifically the reformulation of a VI as a constrained continuously differentiable optimization problem via the regularized gap function \cite{fukushima1992equivalent} and as an unconstrained continuously differentiable optimization problem via the \(D\)-gap function \cite{yamashita1997unconstrained} (see \cite[\S10]{facchinei2003finite} for a survey and \cite{li2007exact}, \cite{patrinos2011global} for applications to constrained optimization and model predictive control of dynamical systems).

		\subsection{Contributions}
			We propose an algorithm that addresses problem \eqref{eq:GenProb} by means of a Newton-like method on the FBE.
Differently from a direct application of the classical Newton method, our approach does not require twice differentiability of the FBE (which would impose additional properties on \(f\) and \(g\)), but merely twice differentiability of \(f\).
This is possible thanks to the introduction of an \emph{approximate generalized Hessian} which only requires access to \(\nabla^2 f\) and to the generalized (Clarke) Jacobian of the proximal mapping of \(g\), as opposed to third-order derivatives and classical Jacobian, respectively.
Moreover, it allows for inexact solutions of linear systems to compute the update direction, which can be done efficiently with a truncated CG method; in particular, no computation of full (generalized) Hessian matrices is necessary, as only (generalized) directional derivatives are needed.
The method is thus particularly appealing when the Clarke Jacobians are sparse and/or well structured, so that the implementation of CG becomes extremely efficient.
Under an error bound condition and a (semi)smoothness assumption at the limit point, which is not required to be isolated, the algorithm exhibits asymptotic \(Q\)-superlinear rates.
For the reader's convenience we collect explicit formulas of the needed Jacobians of the proximal mapping for a wide range of frequently encountered functions, and discuss when they satisfy the needed semismoothness requirements that enable superlinear rates.

		\subsection{Related work}
			This work is a revised version of the unpublished manuscript \cite{patrinos2014forward} and extends ideas proposed in \cite{patrinos2013proximal}, where the FBE is first introduced.
Other FBE-based algorithms are proposed in \cite{stella2017forward,themelis2018forward,stella2017simple}; differently from the truncated-CG type of approximation proposed here, they all employ quasi-Newton directions to mimick second-order information.
The underlying ideas can also be extended to enhance other popular proximal splitting algorithms: the Douglas Rachford splitting (DRS) and the alternating direction method of multipliers (ADMM) \cite{themelis2017douglas}, and for strongly convex problems also the alternating minimization algorithm (AMA) \cite{stella2018newton}.

The algorithm proposed in this paper adopts the recent techniques investigated in \cite{themelis2018forward,stella2017simple} to enhance and greatly simplify the scheme in \cite{patrinos2014forward}.
In particular, \(Q\)-linear and \(Q\)-superlinear rates of convergence are established under an error bound condition, as opposed to uniqueness of the solution.
The proofs of superlinear convergence with an error bound pattern the arguments in \cite{zhou2005superlinear,zhou2006convergence}, although with less conservative requirements.

		\subsection{Organization}
			The work is structured as follows.
In \Cref{sec:Preliminaries} we introduce the adopted notation and list some known facts on generalized differentiability needed in the sequel.
\Cref{sec:ProxAlg} offers an overview on the connections between FBS and the proximal point algorithm, and serves as a prelude to \Cref{sec:FBE} where the forward-backward envelope function is introduced and analyzed.
\Cref{sec:Algorithm} deals with the proposed truncated-Newton algorithm and its convergence analysis.
In \Cref{sec:GenJac} we collect explicit formulas for the generalized Jacobian of the proximal mapping of a rich list of nonsmooth functions, needed for computing the update directions in the proposed algorithm.
Finally, \Cref{sec:Conclusions} draws some conclusions.

	\section{Preliminaries}\label{sec:Preliminaries}

		\subsection{Notation and known facts}
			Our notation is standard and follows that of convex analysis textbooks \cite{bauschke2017convex,bertsekas2015convex,hiriarturruty2004fundamentals,rockafellar1970convex}.
For the sake of clarity we now properly specify the adopted conventions, and briefly recap known definitions and facts in convex analysis.
The interested reader is referred to the above-mentioned textbooks for the details.
	
\textbf{Matrices and vectors.}~
	The \(n\times n\) identity matrix is denoted as \(\I_n\), and the \(\R^n\) vector with all elements equal to \(1\) is as \({\bf 1}_n\); whenever \(n\) is clear from context we simply write \(\I\) or \(\bf 1\), respectively.
	We use the Kronecker symbol \(\delta_{i,j}\) for the \((i,j)\)-th entry of \(\I\).
	Given \(v\in\R^n\), with \(\diag v\) we indicate the \(n\times n\) diagonal matrix whose \(i\)-th diagonal entry is \(v_i\).
	With \(\sym(\R^n)\), \(\sym_+(\R^n)\) and \(\sym_{++}(\R^n)\) we denote respectively the set of symmetric, symmetric positive semidefinite, and symmetric positive definite matrices in \(\R^{n\times n}\).
	
	The minimum and maximum eigenvalues of \(H\in\sym(\R^n)\) are denoted as \(\lambda_{\rm min}(H)\) and \(\lambda_{\rm max}(H)\), respectively.
	For \(Q,R\in\sym(R^n)\) we write \(Q\succeq R\) to indicate that \(Q-R\in\sym_+(\R^n)\), and similarly \(Q\succ R\) indicates that \(Q-R\in\sym_{++}(\R^n)\).
	Any matrix \(Q\in\sym_+(\R^n)\) induces the semi-norm \(\|{}\cdot{}\|_Q\) on \(\R^n\), where \(\|x\|_Q^2\coloneqq\innprod{x}{Qx}\); in case \(Q=\I\), that is, for the Euclidean norm, we omit the subscript and simply write \(\|{}\cdot{}\|\).
	No ambiguity occurs in adopting the same notation for the induced matrix norm, namely \(\|M\|\coloneqq\max\set{\|Mx\|}[x\in\R^n,~\|x\|=1]\) for \(M\in\R^{n\times n}\).

\textbf{Topology.}~
	The \DEF{convex hull} of a set \(E\subseteq\R^n\), denoted as \(\conv E\), is the smallest convex set that contains \(E\) (the intersection of convex sets is still convex).
	The \DEF{affine hull} \(\aff E\) and the \DEF{conic hull} \(\cone E\) are defined accordingly.
	Specifically,
	\begin{align*}
		\conv E
	{}\coloneqq{} &
		\textstyle
		\set{\sum_{i=1}^k\alpha_ix_i}[
			k\in\N,\ 
			x_i\in E,\ 
			\alpha_i\geq0,\ 
			\sum_{i=1}^k\alpha_i=1
		],
	\\
		\cone E
	{}\coloneqq{} &
		\textstyle
		\set{\sum_{i=1}^k\alpha_ix_i}[
			k\in\N,\ 
			x_i\in E,\ 
			\alpha_i\geq 0
		],
	\\
		\aff E
	{}\coloneqq{} &
		\textstyle
		\set{\sum_{i=1}^k\alpha_ix_i}[
			k\in\N,\ 
			x_i\in E,\ 
			\alpha_i\in\R,\ 
			\sum_{i=1}^k\alpha_i=1
		].
	\end{align*}
	The \DEF{closure} and \DEF{interior} of \(E\) are denoted as \(\cl E\) and \(\interior E\), respectively, whereas its \DEF{relative interior}, namely the interior of \(E\) as a subspace of \(\aff E\), is denoted as \(\relint E\).
	With \(\ball xr\) and \(\cball xr\) we indicate, respectively, the open and closed balls centered at \(x\) with radius \(r\).
	
\textbf{Sequences.}~
	The notation \(\seq{a^k}[k\in K]\) represents a sequence indexed by elements of the set \(K\), and given a set \(E\) we write \(\seq{a^k}[k\in K]\subset E\) to indicate that \(a^k\in E\) for all indices \(k\in K\).
	We say that \(\seq{a^k}[k\in K]\subset\R^n\) is \DEF{summable} if \(\sum_{k\in K}\|a^k\|\) is finite, and \DEF{square-summable} if \(\seq{\|a^k\|^2}[k\in K]\) is summable.
	We say that the sequence converges to a point \(a\in\R^n\) \DEF{superlinearly} if either \(a^k=a\) for some \(k\in\N\), or \(\nicefrac{\|a^{k+1}-a\|}{\|a^k-a\|}\to0\); if \(\nicefrac{\|a^{k+1}-a\|}{\|a^k-a\|^q}\) is bounded for some \(q>1\), then we say that the sequence converges \DEF{superlinearly with order \(q\)}, and in case \(q=2\) we say that the convergence is \DEF{quadratic}.

\textbf{Extended-real valued functions.}~
	The extended-real line is \(\Rinf=\R\cup\set{\infty}\).
	Given a function \(\func{h}{\R^n}{[-\infty,\infty]}\), its \DEF{epigraph} is the set
	\begin{align*}
		\epi h
	{}\coloneqq{} &
		\set{(x,\alpha)\in\R^n\times\R}[h(x)\leq\alpha],
	\shortintertext{while its \DEF{domain} is}
		\dom h
	{}\coloneqq{} &
		\set{x\in\R^n}[
			h(x)<\infty
		],
	\shortintertext{and for \(\alpha\in\R\) its \DEF{\(\alpha\)-level set} is}
		\lev_{\leq\alpha}h
	{}\coloneqq{} &
		\set{x\in\R^n}[
			h(x)\leq\alpha
		].
	\end{align*}
	Function \(h\) is said to be \DEF{lower semicontinuous (lsc)} if \(\epi h\) is a closed set in \(\R^{n+1}\) (equivalently, \(h\) is said to be \DEF{closed}); in particular, all level sets of an lsc function are closed.
	We say that \(h\) is \DEF{proper} if \(h>-\infty\) and \(\dom h\neq\emptyset\), and that it is \DEF{level bounded} if for all \(\alpha\in\R\) the level set \(\lev_{\leq\alpha}h\) is a bounded subset of \(\R^n\).

	\textbf{Continuity and smoothness.}~
		A function \(\func{G}{\R^n}{\R^m}\) is \DEF{\(\vartheta\)-H\"older continuous} for some \(\vartheta>0\) if there exists \(L\geq 0\) such that
		\[
			\|G(x)-G(x')\|
		{}\leq{}
			L\|x-x'\|^\vartheta
		\]
		for all \(x,x'\).
		In case \(\vartheta=1\) we say that \(G\) is (\(L\)-)Lipschitz continuous.
		\(G\) is \DEF{strictly differentiable} at \(\bar x\in\R^n\) if the Jacobian matrix
		\(
			JG(\bar x)
		{}\coloneqq{}
			\bigl(\dep{G_i}{x_j}(\bar x)\bigr)_{i,j}
		\)
		exists and
		\[
			\lim_{\substack{x,x'\to\bar x\\x\neq x'}}{
				\frac{
					\|G(x')-JG(\bar x)(x'-x)-G(x)\|
				}{
					\|x'-x\|
				}
			}
		{}={}
			0.
		\]	
	%
		The class of functions \(\func{h}{\R^n}{\R}\) that are \(k\) times continuously differentiable is denoted as \(C^k(\R^n)\).
		We write \(h\in C^{1,1}(\R^n)\) to indicate that \(h\in C^1(\R^n)\) and that \(\nabla h\) is Lipschitz continuous with modulus \(L_h\).
		To simplify the terminology, we will say that such an \(h\) is \DEF{\(L_h\)-smooth}.
		If \(h\) is \(L_h\)-smooth and convex, then for any \(u,v\in\R^n\)
		\begin{equation}\label{eq:LipBound}
			0
		{}\leq{}
			h(v)
			{}-{}
			\bigl[
				h(u)+\innprod{\nabla h(u)}{v-u}
			\bigr]
		{}\leq{}
			\tfrac{L_h}{2}\|v-u\|^2.
		\end{equation}
		Moreover, having \(h\in C^{1,1}(\R^n)\) and \(\mu_h\)-strongly convex is equivalent to having
		\begin{equation}\label{eq:Lmu}
			\mu_h\|v-u\|^2
		{}\leq{}
			\innprod{\nabla h(v)-\nabla h(u)}{v-u}
		{}\leq{}
			L_h\|v-u\|^2
		\end{equation}
		for all \(u,v\in\R^n\).

\textbf{Set-valued mappings.}~
	We use the notation \(\ffunc H{\R^n}{\R^m}\) to indicate a point-to-set function \(\func H{\R^n}{\mathcal P(\R^m)}\), where \(\mathcal P(\R^m)\) is the power set of \(\R^m\) (the set of all subsets of \(\R^m\)).
	The \DEF{graph} of \(H\) is the set
	\begin{align*}
		\graph H
	{}\coloneqq{} &
		\set{(x,y)\in\R^n\times\R^m}[y\in H(x)],
	\shortintertext{while its \DEF{domain} is}
		\dom H
	{}\coloneqq{} &
		\set{x\in\R^n}[
			H(x)\neq\emptyset
		].
	\end{align*}
	We say that \(H\) is \DEF{outer semicontinuous (osc)} at \(\bar x\in\dom H\) if for any \(\varepsilon>0\) there exists \(\delta>0\) such that
	\(
		H(x)\subseteq H(\bar x)+\ball0\varepsilon
	\)
	for all \(x\in\ball{\bar x}{\delta}\).
	In particular, this implies that whenever \(\seq{x^k}\subseteq\dom H\) converges to \(x\) and \(\seq{y^k}\) converges to \(y\) with \(y^k\in H(x^k)\) for all \(k\), it holds that \(y\in H(x)\).
	We say that \(H\) is osc (without mention of a point) if \(H\) is osc at every point of its domain or, equivalently, if \(\graph H\) is a closed subset of \(\R^n\times\R^m\) (notice that this notion does \emph{not} reduce to lower semicontinuity for a single-valued function \(H\)).

\textbf{Convex analysis.}~
	The \DEF{indicator function} of a set \(S\subseteq\R^n\) is denoted as \(\func{\indicator_S}{\R^n}{\Rinf}\), namely
	\[
		\indicator_S(x)
	{}={}
		\begin{cases}[l @{~~} l]
			0 & \text{if }x\in S,
		\\
			\infty & \text{otherwise.}
		\end{cases}
	\]
	If \(S\) is nonempty closed and convex, then \(\indicator_S\) is proper convex and lsc, and both the projection \(\func{\proj_S}{\R^n}{\R^n}\) and the distance \(\func{\dist({}\cdot{},S)}{\R^n}{[0,\infty)}\) are well-defined functions, given by \(\proj_S(x)=\argmin_{z\in S}\|z-x\|\) and \(\dist(x,S)=\min_{z\in S}\|z-x\|\), respectively.
	
	The \DEF{subdifferential} of \(h\) is the set-valued mapping \(\ffunc{\partial h}{\R^n}{\R^n}\) defined as
	\[
		\partial h(x)
	{}\coloneqq{}
		\set{v\in\R^n}[
			h(z)\geq h(x)
			{}+{}
			\innprod{v}{z-x}
			~~
			\forall z\in\R^n
			\vphantom{\seq{x^k}}
		].
	\]
	A vector \(v\in\partial h(x)\) is called a \DEF{subgradient} of \(h\) at \(x\).
	It holds that \(\dom\partial h\subseteq\dom h\), and if \(h\) is proper and convex, then \(\dom\partial h\) is a nonempty convex set containing \(\relint\dom h\), and \(\partial h(x)\) is convex and closed for all \(x\in\R^n\).
	
	A function \(h\) is said to be \(\mu\)-strongly convex for some \(\mu\geq 0\) if \(h-\tfrac\mu2\|{}\cdot{}\|^2\) is convex.
	Unless differently specified, we allow for \(\mu=0\) which corresponds to \(h\) being convex but not strongly so.
	If \(\mu>0\), then \(h\) has a unique (global) minimizer.

		\subsection{Generalized differentiability}
			Due to its inherent nonsmooth nature, classical notions of differentiability may not be directly applicable in problem \eqref{eq:GenProb}.
This subsection contains some definitions and known facts on generalized differentiability that will be needed later on in the paper.
The interested reader is referred to the textbooks \cite{clarke1990optimization,facchinei2003finite,rockafellar2011variational} for the details.
\begin{defin}[Bouligand and Clarke subdifferentials]\label{def:Jacs}%
	Let \(\func{G}{\R^n}{\R^m}\) be locally Lipschitz continuous, and let \(C_G\subseteq\R^n\) be the set of points at which \(G\) is differentiable (in particular \(\R^n\setminus C_G\) has measure zero).
	The \(B\)-subdifferential (also known as \DEF{Bouligand} or \DEF{limiting Jacobian}) of \(G\) at \(\bar x\) is the set-valued mapping
	\(
		\ffunc{\partial_BG}{\R^n}{\R^{m\times n}}
	\)
	defined as
	\[
		\partial_BG(\bar x)
	{}\coloneqq{}
		\set{H\in\R^{m\times n}}[
			\exists\seq{x^k}\subset C_G
			\text{ with }
			x^k\to\bar x,
			JG(x^k)\to H
		],
	\]
	whereas the (\DEF{Clarke}) \DEF{generalized Jacobian} of \(G\) at \(\bar x\) is
	\(
		\ffunc{\partial_CG}{\R^n}{\R^{m\times n}}
	\)
	given by
	\[
		\partial_CG(\bar x)\coloneqq\conv(\partial_BG(\bar x)).
	\]
\end{defin}
If \(\func G{\R^n}{\R^m}\) is locally Lipschitz on \(\R^n\), then \(\partial_CG(x)\) is a nonempty, convex and compact subset of \(\R^{m\times n}\) matrices, and as a set-valued mapping it is osc at every \(x\in\R^n\).
\emph{Semismooth} functions \cite{qi1993nonsmooth} are precisely Lipschitz-continuous mappings for which the generalized Jacobian (and consequenlty the \(B\)-subdifferential) furnishes a first-order  approximation.
\begin{defin}[Semismooth mappings]\label{def:Semismooth}%
	Let \(\func G{\R^n}{\R^m}\) be locally Lipschitz continuous  at \(\bar x\).
	We say that \(G\) is \DEF{semismooth} at \(\bar x\) if
	\begin{subequations}
		\begin{align}
		\label{eq:Semismooth}
			\limsup_{
				\substack{
					x\to\bar x\\
					H\in\partial_CG(x)
				}
			}{
				\frac{
					\|G(x)+H(\bar x-x)-G(\bar x)\|
				}{
					\|x-\bar x\|
				}
			}
		{}={} &
			0.
		\shortintertext{%
			We say that \(G\) is \DEF{\(\vartheta\)-order semismooth} for some \(\vartheta>0\) if the condition can be strengthened to
		}
		\label{eq:thetaSemismooth}
			\limsup_{
				\substack{
					x\to\bar x\\
					H\in\partial_CG(x)
				}
			}{
				\frac{
					\|G(x)+H(\bar x-x)-G(\bar x)\|
				}{
					\|x-\bar x\|^{1+\vartheta}
				}
			}
		{}<{} &
			\infty,
		\end{align}
	and in case \(\vartheta=1\) we say that \(G\) is \DEF{strongly semismooth}.
	\end{subequations}
\end{defin}
To simplify the notation, we adopt the small-\(o\) and big-\(O\) convention to write expressions as \eqref{eq:Semismooth} in the compact form
\(
	G(x)+H(\bar x-x)-G(\bar x)
{}={}
	o(\|x-\bar x\|)
\),
and similarly \eqref{eq:thetaSemismooth} as
\(
	G(x)+H(\bar x-x)-G(\bar x)
{}={}
	O(\|x-\bar x\|^{1+\vartheta})
\).
We remark that the original definition of semismoothness given by \cite{mifflin1977semismooth} requires \(G\) to be directionally differentiable at \(x\).
The definition given here is the one employed by \cite{gowda2004inverse}.
It is also worth remarking that \(\partial_C G(x)\) can be replaced with the smaller set \(\partial_B G(x)\) in \Cref{def:Semismooth}.
Fortunately, the class of semismooth mappings is rich enough to include many functions arising in interesting applications.
For example \emph{piecewise smooth (\(PC^1\)) mappings} are semismooth everywhere.
Recall that a continuous mapping \(\func G{\R^n}{\R^m}\) is \(PC^1\) if there exists a finite collection of smooth mappings \(\func{G_i}{\R^n}{\R^m}\), \(i=1,\ldots,N\), such that
\[
	G(x)\in\set{G_1(x),\ldots,G_N(x)}
\quad
	\forall x\in\R^n.
\]
The definition of \(PC^1\) mapping given here is less general than the one of, e.g., \cite[\S4]{scholtes2012piecewise} but it suffices for our purposes.
For every \(x\in\R^n\) we introduce the set of essentially active indices
\[
	I_G^e(x)
{}\coloneqq{}
	\set{i
	}[{
		x\in\cl\bigl(
			\interior\set{w}[G(w)=G_i(w)]
		\bigr)
	}].
\]
In other words, \(I_G^e(x)\) contains only indices of the pieces \(G_i\) for which there exists a full-dimensional set on which \(G\) agrees with \(G_i\).
In accordance to \Cref{def:Jacs}, the generalized Jacobian of \(G\) at \(x\) is the convex hull of the Jacobians of the essentially active pieces, \ie \cite[Prop. 4.3.1]{scholtes2012piecewise}
\begin{equation}\label{eq:PCJac}
	\partial_C G(x)
{}={}
	\conv\set{JG_i(x)}[i\in I_G^e(x)].
\end{equation}
The following definition is taken from \cite[Def. 7.5.13]{facchinei2003finite}.
\begin{defin}[Linear Newton approximation]\label{def:LNA}%
	Let \(\func G{\R^n}{\R^m}\) be continuous on \(\R^n\).
	We say that \(G\) admits a \DEF{linear Newton approximation (LNA)} at \(\bar x\in\R^n\) if there exists a set-valued mapping \(\ffunc{\mathcal H}{\R^n}{\R^{m\times n}}\) that has nonempty compact images, is outer semicontinuous at \(\bar x\), and
	\begin{align*}
		\limsup_{
			\substack{
				\mathllap x{}\to{}\mathrlap{\bar x}\\
				\mathllap H\fillwidthof[c]{{}\to{}}{\in}\mathrlap{\mathcal H(x)}
			}
		}{
			\frac{
				\|G(x)+H(\bar x-x)-G(\bar{x})\|
			}{
				\|x-\bar x\|
			}
		}
	{}={} &
		0.
	\shortintertext{%
		If for some \(\vartheta>0\) the condition can be strengthened to
	}
		\limsup_{
			\substack{
				\mathllap x{}\to{}\mathrlap{\bar x}\\
				\mathllap H\fillwidthof[c]{{}\to{}}{\in}\mathrlap{\mathcal H(x)}
			}
		}{
			\frac{
				\|G(x)+H(\bar x-x)-G(\bar{x})\|
			}{
				\|x-\bar x\|^{1+\vartheta}
			}
		}
	{}<{} &
		\infty,
	\end{align*}
	then we say that \(\mathcal H\) is a \DEF{\(\vartheta\)-order LNA}, and if \(\vartheta=1\) we say that \(\mathcal H\) is a \DEF{strong LNA}.
\end{defin}
Functions \(G\) as in \Cref{def:LNA} are also referred to as \(\mathcal H\)-semismooth in the literature, see \eg \cite{ulbrich2009optimization}, however we prefer to stick to the terminology of \cite{facchinei2003finite} and rather say that \(\mathcal H\) is a LNA for \(G\).
Arguably the most notable example of a LNA for semismooth mappings is the generalized Jacobian, cf. \Cref{def:Jacs}.
However, semismooth mappings can admit LNAs different from the generalized Jacobian.
More importantly, mappings that are not semismooth may also admit a LNA.

\begin{lem}[{\cite[Prop. 7.4.10]{facchinei2003finite}}]\label{thm:fp_prop7.4.10}%
	Let \(h\in C^1(\R^n)\) and suppose that \(\ffunc{\mathcal H}{\R^n}{\R^{n\times n}}\) is a LNA for \(\nabla h\) at \(\bar x\).
	Then,

	\[
		\lim_{
			\substack{
				\mathllap x{}\to{}\mathrlap{\bar x}\\
				\mathllap H\fillwidthof[c]{{}\to{}}{\in}\mathrlap{\mathcal H(x)}
			}
		}{
			\frac{
				h(x)-h(\bar x)
				{}-{}
				\innprod{\nabla h(\bar x)}{x-\bar x}
				{}-{}
				\tfrac12\innprod{H(x-\bar x)}{x-\bar x}
			}{
				\|x-\bar x\|^2
			}
		}
	{}={}
		0.
	\]
\end{lem}
We remark that although \cite[Prop. 7.4.10]{facchinei2003finite} assumes semismoothness of \(\nabla h\) at \(\bar x\) and uses \(\partial_C(\nabla h)\) in place of \(\mathcal H\); however, exactly the same arguments apply for any LNA of \(\nabla h\) at \(\bar x\) even without the semismoothness assumption.

	\section{Proximal algorithms}\label{sec:ProxAlg}%

		\subsection{Proximal point and Moreau envelope}\label{sec:Prox}%
\begin{figure}
		\begin{minipage}[t][][t]{0.45\linewidth}
			\vspace*{-\abovecaptionskip}%
			\caption{%
				\emph{%
					Moreau envelope of the function\newline
					\hspace*{0pt}\hfill\(
						\varphi(x) = \tfrac13x^3+x^2-x+1+\indicator_{[0,\infty)}(x)
					\)\hfill\hspace*{0pt}\newline
					with parameter \(\gamma=0.2\).%
				}
				\newline
				At each point \(x\), the Moreau envelope \(\varphi^\gamma\) is the minimum of the quadratic majorization model
				\(
					\mathcal M_\gamma^\varphi=\varphi+\frac{1}{2\gamma}({}\cdot{}-x)^2
				\),
				the unique minimizer being, by definition, the proximal point \(\hat x\coloneqq\prox_{\gamma\varphi}(x)\).
				It is a convex smooth lower bound to \(\varphi\), despite the fact that \(\varphi\) might be extended-real valued.
				Function \(\varphi\) and its Moreau envelope \(\varphi^\gamma\) have same \(\inf\) and \(\argmin\); in fact, the two functions agree (only) on the set of minimizers.
				In general, \(\varphi^\gamma\) is \emph{sandwiched} as
				\(
					\varphi\circ\prox_{\gamma \varphi}\leq\varphi^\gamma\leq\varphi
				\).
			}%
			\label{fig:ME}%
		\end{minipage}
		\hfill
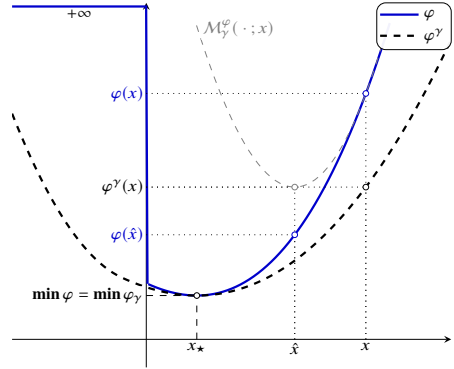
		\begin{minipage}[t][][t]{0.5\linewidth}
			\vspace*{0pt}%
			{{%
			\pgfkeys{/pgf/images/include external/.code={\includegraphics[width=\linewidth]{#width=\linewidth}}}%
			\tikzsetnextfilename{PP}%
			\begin{tikzpicture}
\def\gam{0.2}
	\def\xMin{-1.1}
	\def\xMax{2.5}
	\def\yMin{-0.5}
	\def\yMax{6}
	\def\xo{1.8}
	\def\a{1}
	\def\b{-1}
	\def\c{1}
	\def\xstar{max(-\a+sqrt(\a^2-\b), 0)}
	\tikzset{%
		declare function = {%
			phi(\x)     = (\x < 0) * \YINFTY  +  (\x >= 0) * min(\x^3/3 + \a*\x^2 + \b*\x + \c, \YINFTY);
			proxphiunconstrained(\x) = (-(2*\gam*\a + 1) + sqrt( (2*\gam*\a + 1)^2 - 4*(\b*\gam-\x)*\gam ) ) / (2*\gam);
			proxphi(\x) = (proxphiunconstrained(\x) >= 0) * proxphiunconstrained(\x);
		},
	}
	\tikzset{%
		declare function = {
			phiME(\x) = phi(proxphi(\x)) + (\x-proxphi(\x))^2 / (2*\gam);
			Q(\x,\w) = (phi(\x) + (\x-\w)^2/(2*\gam));
		},
	}
	\begin{axis}[myaxis]
		%
		\plotExtFunc[phi]{phi}{\xMin:0.01}
		\addlegendentry{\(\varphi\)}
		\node[label, anchor=north] at (axis cs: {(\xMin)/2},{\YINFTY}) {\(+\infty\)};
		%
		\addplot[Q, dashed, forget plot] {Q(x,\xo)};
		\node[label, color=Q, anchor=south west, outer sep=0pt] at (axis cs: {0.25*\xo},{Q(0.25*\xo,\xo)}) {\(\mathcal M_\gamma^\varphi({}\cdot{};x)\)};
		%
		\addplot[ME] {phiME(x)};
		\addlegendentry{\(\varphi^\gamma\)}
		%
		\node[nosep] (phi0) at (axis cs: {\xo},{phi(\xo)}) {};
		\drawCoord{phi0}
		\drawDot[color=phi]{phi0};
		\xLabel{phi0}{\(x\vphantom{x_\star}\)}
		\yLabel[color=phi]{phi0}{$\varphi(x)$}
		%
		\node[nosep] (me0) at (axis cs: {\xo},{phiME(\xo)}) {};
		\drawCoord{me0}
		\drawDot[color=ME]{me0};
		\yLabel[color=ME]{me0}{\(\varphi^\gamma(x)\)}
		%
		\node[nosep] (Q0) at (axis cs: {proxphi(\xo)},{phiME(\xo)}) {};
		\drawCoord{Q0}
		\drawDot[color=Q]{Q0}
		\xLabel{Q0}{\(\hat x\vphantom{x_\star}\)}
		%
		\node[nosep] (phiP0) at (axis cs: {proxphi(\xo)},{phi(proxphi(\xo))}) {};
		\drawCoord{phiP0}
		\drawDot[color=phi]{phiP0}
		\yLabel[color=phi]{phiP0}{\(\varphi(\hat x)\)}
		%
		\node[nosep] (xmin) at (axis cs: {\xstar},{phi(\xstar)}) {};
		\drawCoord[dashed][dashed]{xmin}
		\drawDot{xmin}
		\xLabel[anchor=south]{xmin}{\(x_\star\vphantom{x_\star}\)}
		\yLabel[anchor=north east]{xmin}{\(\min\varphi=\min\FBE\)}
	\end{axis}
\end{tikzpicture}%
		}}%
		\end{minipage}
	\end{figure}

				The \DEF{proximal mapping} of a proper closed and convex function \(\func h{\R^n}{\Rinf}\) with parameter \(\gamma>0\) is \(\func{\prox_{\gamma h}}{\R^n}{\R^n}\), given by
	\begin{align}\label{eq:prox}
		\prox_{\gamma h}(x)
	{}\coloneqq{} &
		\argmin_{w\in\R^n}\set{
			\smashoverbracket{
				h(w)+\tfrac{1}{2\gamma}\|w-x\|^2
			}{
				\mathcal M_\gamma^h(w;x)
			}
		}.
	\shortintertext{%
		The \emph{majorization model} \(\mathcal M_\gamma^h(x;{}\cdot{})\) is a proper and strongly convex function, and therefore has a unique minimizer.
		The value function, as opposed to the minimizer, defines the \DEF{Moreau envelope} \(\func{h^\gamma}{\R^n}{\R}\), namely%
	}
	\label{eq:ME}
		h^\gamma(x)
	{}\coloneqq{} &
		\min_{w\in\R^n}\set{
			h(w)+\tfrac{1}{2\gamma}\|w-x\|^2
		},
	\end{align}
	which is real valued and Lipschitz differentiable, despite the fact that \(h\) might be extended-real valued.
	Properties of the Moreau envelope and the proximal mapping are well documented in the literature, see \eg \cite[\S24]{bauschke2017convex}.
	For example, \(\prox_{\gamma h}\) is nonexpansive (Lipschitz continuous with modulus \(1\)) and is characterized by the implicit inclusion
	\begin{equation}\label{eq:proxCharacterization}
		\hat x=\prox_{\gamma h}(x)
	\quad\Leftrightarrow\quad
		\tfrac1\gamma(x-\hat x)\in\partial h(\hat x).
	\end{equation}
	For the sake of a brief recap, we now list some other important known properties.
	\Cref{thm:MEsandwich} provides some relations between \(h\) and its Moreau envelope \(h^\gamma\), which we informally refer to as \emph{sandwich property} for apparent reasons, cf. \Cref{fig:ME}.
	\Cref{thm:MEequiv} highlights that the minimization of a (proper, lsc and) convex function can be expressed as the convex smooth minimization of its Moreau envelope.
	\begin{thm}[Moreau envelope: sandwich property \cite{bauschke2017convex,chen1993convergence}]\label{thm:MEsandwich}%
		For all \(\gamma>0\) the following hold for the cost function \(\varphi\):
		\begin{enumerate}
		\item\label{thm:MEleq}%
			\(
				\varphi^\gamma(x)
			{}\leq{}
				\fillwidthof[c]{\varphi^\gamma(x)}{\varphi(x)}
				{}-{}
				\frac{1}{2\gamma}\|x-\hat x\|^2
			\)
			~for all~
			\(x\in\R^n\)
			~where~
			\(\hat x\coloneqq\prox_{\gamma\varphi}(x)\);
		\item\label{thm:MEgeq}%
			\(
				\hphantom{\varphi^\gamma(x)}
				\mathllap{\varphi(\hat x)}
			\fillwidthof[c]{{}\leq{}}{{}={}}
				\varphi^\gamma(x)
				{}-{}
				\frac{1}{2\gamma}\|x-\hat x\|^2
			\)
			~for all~
			\(x\in\R^n\)
			~where~
			\(\hat x\coloneqq\prox_{\gamma\varphi}(x)\);
		\item\label{thm:ME=}%
			\(
				\varphi^\gamma(x)
			\fillwidthof[c]{{}\leq{}}{{}={}}
				\fillwidthof[c]{\varphi^\gamma(x)}{\varphi(x)}
			\)
			~iff~
			\(
				x\in\argmin\varphi
			\).
		\end{enumerate}
		\begin{proof}
			\begin{proofitemize}
			\item\ref{thm:MEleq}.
				This fact is shown in \cite[Lem. 3.2]{chen1993convergence} for a more general notion of proximal point operator; namely, the square Euclidean norm appearing in \eqref{eq:prox} and \eqref{eq:ME} can be replaced by arbitrary Bregman divergences.
				In this simpler case, since \(\tfrac1\gamma(x-\hat x)\) is a subgradient of \(\varphi\) at \(\hat x\), cf. \eqref{eq:proxCharacterization}, we have
				\[
					\varphi(x)
				{}\geq{}
					\varphi(\hat x)
					{}+{}
					\innprod{\tfrac1\gamma(x-\hat x)}{x-\hat x}
				{}={}
					\varphi(\hat x)
					{}+{}
					\tfrac1\gamma\|x-\hat x\|^2.
				\]
				The claim now follows by substracting \(\tfrac{1}{2\gamma}\|x-\hat x\|^2\) from both sides.
			\item\ref{thm:MEgeq}.
				Follows by definition, cf. \eqref{eq:prox} and \eqref{eq:ME}.
			\item\ref{thm:ME=}.
				See \cite[Prop. 17.5]{bauschke2017convex}.
			\qedhere
			\end{proofitemize}
		\end{proof}
	\end{thm}
	\begin{thm}[Moreau envelope: convex smooth minimization equivalence \cite{bauschke2017convex}]\label{thm:MEequiv}%
		For all \(\gamma>0\) the following hold for the cost function \(\varphi\):
		\begin{enumerate}
		\item\label{thm:ME:C1}%
			\(\varphi^\gamma\)
			~is convex and smooth with~
			\(L_{\varphi^\gamma}=\gamma^{-1}\)
			~and~
			\(
				\nabla\varphi^\gamma(x)
			{}={}
				\gamma^{-1}\bigl(x-\prox_{\gamma\varphi}(x)\bigr)
			\);
		\item\label{thm:ME:infequiv}%
			\(
				\inf\varphi=\inf\varphi^\gamma
			\);
		\item\label{thm:ME:minimequiv}%
			\(
				x_\star\in\argmin\varphi
			\)
			~iff~
			\(
				x_\star\in\argmin\varphi^\gamma
			\)
			~iff~
			\(
				\nabla\varphi^\gamma(x_\star)=0
			\).
		\end{enumerate}
		\begin{proof}
			\begin{proofitemize}
			\item\ref{thm:ME:C1}.
				See \cite[Prop.s 12.15 and 12.30]{bauschke2017convex}.
			\item\ref{thm:ME:infequiv}.
				See \cite[Prop. 12.9(iii)]{bauschke2017convex}.
			\item\ref{thm:ME:minimequiv}.
				See \cite[Prop. 17.5]{bauschke2017convex}.
			\qedhere
			\end{proofitemize}
		\end{proof}
	\end{thm}
	As a consequence of \Cref{thm:MEequiv}, one can address the minimization of the convex but possibly nonsmooth and extended-real-valued function \(\varphi\) by means of gradient descent on the smooth envelope function \(\varphi^\gamma\) with stepsize \(0<\tau<\nicefrac{2}{L_{\FBE}}=2\gamma\).
	As first noticed by Rockafellar \cite{rockafellar1976monotone}, this simply amounts to (relaxed) fixed-point iterations of the proximal point operator, namely
	\begin{equation}\label{eq:PPA}
		x^+=(1-\lambda)x+\lambda\prox_{\gamma\varphi}(x),
	\end{equation}
	where \(\lambda=\nicefrac\tau\gamma\in(0,2)\) is a possible relaxation parameter.
	The scheme, known as \emph{proximal point algorithm} (PPA) and first introduced by Martinet \cite{martinet1970breve}, is well covered by the broad theory of monotone operators, where convergence properties can be easily derived with simple tools of Fejérian monotonicity, see \eg \cite[Thm.s 23.41 and 27.1]{bauschke2017convex}.
	Nevertheless, not only does the interpretation as gradient method provide a beautiful theoretical link, but it also enables the employment of acceleration techniques exclusively stemming from smooth unconstrained optimization, such as Nesterov's extrapolation \cite{guler1992newproximal} or quasi-Newton schemes \cite{chen1999proximal}, see also \cite{bertsekas1982constrained} for extensions to the dual formulation.

		\subsection{Forward-backward splitting}
			While it is true that every convex minimization problem can be smoothened by means of the Moreau envelope, unfortunately it is often the case that the computation of the proximal operator (which is needed to evaluate the envelope) is as hard as solving the original problem.
For instance, evaluating the Moreau envelope of the cost of modeling a convex QP at one point amounts to solving another QP with same constraints and augmented cost.
To overcome this limitation there comes the idea of \emph{splitting schemes}, which decompose a complex problem in small components which are easier to operate onto.
A popular such scheme is the \DEF{forward-backward splitting} (FBS), which addresses minimization problems of the form \eqref{eq:GenProb}.

Given a point \(x\in\R^n\), one iteration of \DEF{forward-backward splitting} (FBS) for problem \eqref{eq:GenProb} with stepsize \(\gamma>0\) and relaxation \(\lambda>0\) consists in
\begin{align}\label{eq:FBS}
	x^+
{}={} &
	(1-\lambda)x+\lambda \T(x),
\shortintertext{%
	where%
}
\label{eq:PG}
	\T(x)
{}\coloneqq{} &
	\FB x
\end{align}
is the \DEF{forward-backward operator}, characterized as
\begin{equation}\label{eq:FBCharacterization}
	\bar x=T_\gamma(x)
\quad\Leftrightarrow\quad
	\tfrac1\gamma(x-\bar x)-(\nabla f(x)-\nabla f(\bar x))
	{}\in{}
	\partial\varphi(\bar x),
\end{equation}
as it follows from \eqref{eq:proxCharacterization}.
FBS interleaves a gradient descent step on \(f\) and a proximal point step on \(g\), and as such it is also known as \DEF{proximal gradient method}.
If both \(f\) and \(g\) are (lsc, proper and) convex, then the solutions to \eqref{eq:GenProb} are exactly the fixed points of the forward-backward operator \(T_\gamma\).
In other words,
\begin{equation}\label{eq:R0}
	x_\star\in\argmin\varphi
\quad\text{iff}\quad
	\Res(x_\star)=0,
\end{equation}
where
\begin{equation}\label{eq:R}
	\Res(x)
{}\coloneqq{}
	\tfrac1\gamma\bigl(
		x-\FB x
	\bigr)
\end{equation}
is the \DEF{forward-backward residual}.\footnote{%
	Due to apparent similarities with gradient descent iterations, having \(x^+=x-\gamma\Res(x)\) in FBS, \(\Res\) is also referred to as (generalized) \DEF{gradient mapping}, see \eg \cite{drusvyatskiy2018error}.
	In particular, if \(g=0\) then \(\Res=\nabla f\) whereas if \(f=0\) then \(\Res=\nabla g^\gamma\).
	The analogy will be supported by further evidence in the next section where we will see that, up to a change of metric, indeed \(\Res\) is the gradient of the \emph{forward-backward envelope function}.
}
FBS iterations \eqref{eq:FBS} are well known to converge to a solution to \eqref{eq:GenProb} provided that \(f\) is smooth and that the parameters are chosen as \(\gamma\in(0,\nicefrac{2}{L_f})\) and \(\lambda\in(0,2-\nicefrac{\gamma L_f}{2})\) \cite[Cor. 28.9]{bauschke2017convex} (\(\lambda=1\), which is always feasible, is the typical choice).

			\subsection{Error bounds and quadratic growth}
				We conclude the section with some inequalities that will be useful in the sequel.
\begin{lem}\label{thm:cost_dd}%
	Suppose that \(\mathcal X_\star\) is nonempty.
	Then,
	\[
		\varphi(x)-\varphi_\star
	{}\leq{}
		\dist(0,\partial\varphi(x))
		\dist(x,\mathcal X_\star)
	\quad
		\forall x\in\R^n.
	\]
	\begin{proof}
		From the subgradient inequality it follows that for all \(x_\star\in\mathcal X_\star\) and \(v\in\partial\varphi(x)\) we have
		\[
			\varphi(x)-\varphi_\star
		{}={}
			\varphi(x)-\varphi(x_\star)
		{}\leq{}
			\innprod{v}{x-x_\star}
		{}\leq{}
			\|v\|\|x-x_\star\|
		\]
		and the claimed inequality follows from the arbitrarity of \(x_\star\) and \(v\).
	\end{proof}
\end{lem}
\begin{lem}\label{thm:r_d}%
	Suppose that \(\mathcal X_\star\) is nonempty.
	For all \(x\in\R^n\) and \(\gamma>0\) the following holds
	\[
		\|\Res(x)\|
	{}\geq{}
		\tfrac{1}{1+\gamma L_f}
		\dist\bigl(0,\partial\varphi(\T(x))\bigr)
	\]
	\begin{proof}
		Let \(\bar x\coloneqq\T(x)\).
		The characterization \eqref{eq:FBCharacterization} of \(\T\) implies that
		\[
			\|\Res(x)\|
		{}\geq{}
			\dist\bigl(0,\partial\varphi(\bar x)\bigr)
			{}-{}
			\|\nabla f(x)-\nabla f(\bar x)\|
		{}\geq{}
			\dist\bigl(0,\partial\varphi(\bar x)\bigr)
			{}-{}
			\gamma L_f\|\Res(x)\|.
		\]
		After trivial rearrangements the sought inequality follows.
	\end{proof}
\end{lem}
Furhter interesting inequalities can be derived if the cost function \(\varphi\) satisfies an \emph{error bound}, which can be regarded as a generalization of strong convexity that does not require uniqueness of the minimizer.
The interested reader is referred to \cite{luo1993error,pang1997error,bauschke2015linear,drusvyatskiy2018error} and references therein for extensive discussions.
\begin{defin}[Quadratic growth and error bound]
	Suppose that \(\mathcal X_\star\neq\emptyset\).
	Given \(\mu,\nu>0\), we say that
	\begin{enumerateq}
	\item
		\(\varphi\) satisfies the \DEF{quadratic growth} with constants \((\mu,\nu)\) if
		\begin{equation}\label{eq:QG}
			\varphi(x)-\varphi_\star
		{}\geq{}
			\tfrac\mu2\dist(x,\mathcal X_\star)^2
		\quad
			\forall x
			{}\in{}
			\lev_{\leq\varphi_\star+\nu}\varphi;
		\end{equation}
	\item
		\(\varphi\) satisfies the \DEF{error bound} with constants \((\mu,\nu)\) if
		\begin{equation}\label{eq:EB}
			\dist(0,\partial\varphi(x))
		{}\geq{}
			\tfrac\mu2\dist(x,\mathcal X_\star)
		\quad
			\forall x
			{}\in{}
				\lev_{\leq\varphi_\star+\nu}\varphi.
		\end{equation}
	\end{enumerateq}
	In case \(\nu=\infty\) we say that the properties are satisfied \emph{globally}.
\end{defin}
\begin{thm}[{\cite[Thm. 3.3]{drusvyatskiy2018error}}]\label{thm:dru}%
	For a proper convex and lsc function, the quadratic growth with constants \((\mu,\nu)\) is equivalent to the error bound with same constants.
\end{thm}
\begin{lem}[Globality of quadratic growth]\label{thm:QGglobal}%
	Suppose that \(\varphi\) satisfies the quadratic growth with constants \((\mu,\nu)\).
	Then, for every \(\nu'>\nu\) it satisfies the quadratic growth with constants \((\mu',\nu')\), where
	\[
		\mu'
	{}\coloneqq{}
		\tfrac\mu2
		\min\set{
			1,\,
			\tfrac{\nu}{\nu'-\nu}
		}.
	\]
	\begin{proof}
		Let \(\nu'>\nu\) be fixed, and let \(x\in\lev_{\leq\varphi_\star+\nu'}\) be arbitrary.
		Since \(\mu'\leq\mu\), the claim is trivial if \(\varphi(x)\leq\varphi_\star+\nu\); we may thus suppose that \(\varphi(x)>\varphi_\star+\nu\).
		Let \(z\) be the projection of \(x\) onto the (nonempty closed and convex) level set \(\lev_{\leq\varphi_\star+\nu}\), and observe that \(\varphi(z)=\varphi_\star+\nu\).
		With \cref{thm:cost_dd,thm:dru} we can upper bound \(\nu\) as
		\begin{equation}\label{eq:nud2}
			\nu
		{}={}
			\varphi(z)-\varphi_\star
		{}\leq{}
			\dist(0,\partial\varphi(z))
			\dist(z,\mathcal X_\star)
		{}\leq{}
			\tfrac2\mu
			\dist(0,\partial\varphi(z))^2.
		\end{equation}
		Moreover, it follows from \cite[Thm. 1.3.5]{hiriarturruty2004fundamentals} that there exists a subgradient \(v\in\partial\varphi(z)\) such that
		\(
			\innprod{v}{x-z}
		{}={}
			\|v\|\|x-z\|
		\).
		Then,
		\begin{align*}
			\varphi(x)
		{}\geq{} &
			\varphi(z)+\innprod{v}{x-z}
		{}={}
			\varphi(z)+\|v\|\|x-z\|
		{}\geq{}
			\varphi(z)+\dist(0,\partial\varphi(z))\|x-z\|
		\\
		\numberthis\label{eq:x-z}
		{}\overrel[\geq]{\eqref{eq:nud2}}{} &
			\varphi(z)
			{}+{}
			\sqrt{\tfrac{\mu\nu}{2}}\|x-z\|.
		\end{align*}
		By substracting \(\varphi(z)\) from the first and last terms we obtain
		\[
			\|x-z\|
		{}\leq{}
			\sqrt{\tfrac{2}{\mu\nu}}
			\bigl(
				\varphi(x)-\varphi(z)
			\bigr)
		{}\leq{}
			\sqrt{\tfrac{2}{\mu\nu}}
			(\nu'-\nu),
		\]
		which implies
		\begin{equation}\label{eq:linquad}
			\|x-z\|
		{}\geq{}
			\sqrt{\tfrac{\mu\nu}{2}}
			\tfrac{1}{\nu'-\nu}
			\|x-z\|^2.
		\end{equation}
		Thus,
		\begin{align*}
			\varphi(x)
			{}-{}
			\varphi_\star
		{}\overrel*[\geq]{\eqref{eq:x-z}}{} &
			\varphi(z)
			{}-{}
			\varphi_\star
			{}+{}
			\sqrt{\tfrac{\mu\nu}{2}}\|x-z\|
		\shortintertext{%
			using the quadratic growth at \(z\) and the inequality \eqref{eq:linquad}
		}
		{}\geq{} &
			\tfrac\mu2\dist(z,\mathcal X_\star)^2
			{}+{}
			\tfrac{\mu\nu}{2(\nu'-\nu)}
			\|x-z\|^2
		\\
		{}\geq{} &
			\tfrac\mu2
			\min\set{
				1,\,\tfrac{\nu}{\nu'-\nu}
			}
			\left[
				\dist(z,\mathcal X_\star)^2
				{}+{}
				\|x-z\|^2
			\right].
		\end{align*}
		By using the fact that \(a^2+b^2\geq\frac12(a+b)^2\) for any \(a,b\in\R\) together with the triangular inequality \(\dist(x,\mathcal X_\star)\leq\|x-z\|+\dist(z,\mathcal X_\star)\), we conclude that
		\(
			\varphi(x)
			{}-{}
			\varphi_\star
		{}\geq{}
			\tfrac{\mu'}{2}
			\dist(x,\mathcal X_\star)^2
		\),
		with \(\mu'\) as in the statement.
		Since \(\mu'\) depends only on \(\mu\), \(\nu\), and \(\nu'\), from the arbitrarity of \(x\in\lev_{\leq\varphi_\star+\nu'}\) the claim follows.
	\end{proof}
\end{lem}

\begin{thm}[{\cite[Cor. 3.6]{drusvyatskiy2018error}}]\label{thm:EB:Res}%
	Suppose that \(\varphi\) satisfies the quadratic growth with constants \((\mu,\nu)\).
	Then, for all \(\gamma\in(0,\nicefrac{1}{L_f})\) and \(x\in\lev_{\leq\varphi_\star+\nu}\varphi\) we have
	\[
		\dist(x,\mathcal X_\star)
	{}\leq{}
		(\gamma+\nicefrac2\mu)(1+\gamma L_f)
		\|\Res(x)\|.
	\]
\end{thm}

	\section{Forward-backward envelope}\label{sec:FBE}
		There are clearly infinte ways of representing the (proper, lsc and) convex function \(\varphi\) in \eqref{eq:GenProb} as the sum of two convex functions \(f\) and \(g\) with \(f\) smooth, and each of these choices leads to a different FBS operator \(T_\gamma\).
If \(f=0\), for instance, then \(T_\gamma\) reduces to \(\prox_{\gamma\varphi}\), and consequently FBS \eqref{eq:FBS} to the PPA \eqref{eq:PPA}.
A natural question then arises, whether a function exists that serves as ``envelope'' for FBS in the same way that \(\FBE\) does for \(\prox_{\gamma\varphi}\).
We will now provide a positive answer to this question by reformulating the nonsmooth problem \eqref{eq:GenProb} as the minimization of a differentiable function.
To this end, the following requirements on \(f\) and \(g\) will be assumed throughout the paper without further mention.
\begin{svgraybox}
	\vspace*{-\baselineskip}
	\begin{ass}[Basic requirements]
		In problem \eqref{eq:GenProb},
		\begin{enumerate}%
		\item
			\(\func{f}{\R^n}{\R}\) is convex, twice continuously differentiable and \(L_f\)-smooth;
		\item
			\(\func{\hphantom{f}\mathllap g}{\R^n}{\Rinf}\) is lsc, proper and convex.
		\end{enumerate}
		\vspace*{-12pt}%
	\end{ass}
\end{svgraybox}
Compared to the classical FBS assumptions, the only additional requirement is twice differentiability of \(f\).
This ensures that the \emph{forward operator} \(x\mapsto\Fw x\) is differentiable; we denote its Jacobian as \(\func{Q_\gamma}{\R^n}{\R^{n\times n}}\), namely
\begin{equation}\label{eq:Q}
	Q_\gamma(x)
{}\coloneqq{}
	\I-\gamma\nabla^2\!f(x).
\end{equation}
Notice that, due to the bound \(\nabla^2\!f(x)\preceq L_f\I\) (which follows from \(L_f\)-smoothness of \(f\), see \cite[Lem. 1.2.2]{nesterov2003introductory}) \(Q_\gamma(x)\) is invertible (in fact, positive definite) whenever \(\gamma<\nicefrac{1}{L_f}\).
Moreover, due to the chain rule and \Cref{thm:ME:C1} we have that
\begin{align*}
	\nabla{}\bigl[g^\gamma\circ(\Fw{})\bigr](x)
{}={} &
	\gamma^{-1}Q_\gamma(x)
	\bigl[\Fw x-\FB x\bigr]
\\
{}={} &
	Q_\gamma(x)\bigl[\Res(x)-\nabla f(x)\bigr].
\end{align*}
Rearranging,
\begin{align*}
	Q_\gamma(x)\Res(x)
{}={} &
	\nabla f(x)-\gamma\nabla^2\!f(x)\nabla f(x)
	{}+{}
	\nabla{}\bigl[g^\gamma\circ(\Fw{})\bigr](x)
\\
{}={} &
	\nabla f(x)
	{}-{}
	\nabla{}\left[\tfrac\gamma2\|\nabla f\|^2\right](x)
	{}+{}
	\nabla{}\bigl[g^\gamma\circ(\Fw{})\bigr](x)
\\
{}={} &
	\nabla{}\left[
		f
		{}-{}
		\tfrac\gamma2\|\nabla f\|^2
		{}+{}
		g^\gamma\circ(\Fw{})
	\right](x)
\end{align*}
we obtain the gradient of a real-valued function, which we define as follows.
\begin{svgraybox}
	\vspace*{-\baselineskip}%
	\begin{defin}[Forward-backward envelope]
		The \DEF{forward-backward envelope} (FBE) for the composite minimization problem \eqref{eq:GenProb} is the function \(\func{\FBE}{\R^n}{\R}\)
		defined as
		\begin{equation}\label{eq:FBE}
			\FBE(x)
		{}\coloneqq{}
			f(x)-\tfrac\gamma2\|\nabla f(x)\|^2+g^\gamma(\Fw x).
			\vspace*{-\belowdisplayskip}
		\end{equation}
	\end{defin}
\end{svgraybox}
\noindent
In the next section we discuss some of the favorable properties enjoyed by the FBE.

		\subsection{Basic properties}\label{sec:Basic}
			We already verified that the FBE is differentiable with gradient
\begin{equation}\label{eq:gradFBE}
	\nabla\FBE(x)
{}={}
	Q_\gamma(x)\Res(x).
\end{equation}
In particular, for \(\gamma<\nicefrac{1}{L_f}\) one obtains that a FBS step is a (scaled) gradient descent step on the FBE, similarly as the relation between Moreau envelope and PPA; namely,
\begin{equation}\label{eq:FBEGD}
	T_\gamma(x)
{}={}
	x
	{}-{}
	\gamma Q_\gamma(x)^{-1}
	\nabla\FBE(x).
\end{equation}

To take the analysis of the FBE one step further, let us consider the equivalent expression of the operator \(T_\gamma\) as
\begin{align}\label{eq:Targmin}
	T_\gamma(x)
{}={} &
	\argmin_{w\in\R^n}\set{
		\smashoverbracket{
			f(x)+\innprod{\nabla f(x)}{w-x}+\tfrac{1}{2\gamma}\|w-x\|^2
			{}+{}
			g(w)
		}{\mathcal M_\gamma^{f,g}(w;x)}
	}.
\shortintertext{%
	Differently from the quadratic model \(\mathcal M_\gamma^\varphi\) in \eqref{eq:prox}, \(\mathcal M_\gamma^{f,g}\) replaces the differentiable component \(f\) with a linear approximation.
	Building upon the idea of the Moreau envelope, instead of the minimizer \(T_\gamma(x)\) we consider the value attained in the subproblem \eqref{eq:Targmin}, and with simple algebra one can easily verify that this gives rise once again to the FBE:
}
	\label{eq:FBEmin}
	\FBE(x)
{}={} &
	\min_{w\in\R^n}\set{
		f(x)+\innprod{\nabla f(x)}{w-x}+\tfrac{1}{2\gamma}\|w-x\|^2
		{}+{}
		g(w)
	}.
\end{align}
Starting from this expression we can easily mirror the properties of the Moreau envelope stated in \Cref{thm:MEsandwich,thm:MEequiv}.
These results appeared in the independent works \cite{nesterov2013gradient} and \cite{patrinos2013proximal}, although the former makes no mention of an ``envelope'' function and simply analyzes the \emph{majorization-minimization} model \(\mathcal M_\gamma^{f,g}\).
\begin{thm}[FBE: sandwich property]\label{thm:sandwich}%
Let \(\gamma>0\) and \(x\in\R^n\) be fixed, and denote \(\bar x=\T(x)\).
The following hold:
	\begin{enumerate}
	\item\label{thm:leq}%
		\(
			\FBE(x)
		{}\leq{}
			\fillwidthof[c]{\FBE(x)}{\varphi(x)}
			{}-{}
			\frac{1}{2\gamma}\|x-\bar x\|^2
		\);
	\item\label{thm:geq}%
		\(
		\FBE(x)
		{}-{}
		\frac{1}{2\gamma}\|x-\bar x\|^2
		{}\leq{}
			\varphi(\bar x)
		{}\leq{}
			\FBE(x)
			{}-{}
			\frac{1-\gamma L_f}{2\gamma}\|x-\bar x\|^2
		\).
	\end{enumerate}
	In particular,
	\begin{enumerate}[resume]
	\item\label{thm:=}%
		\(
			\FBE(x_\star)
		\fillwidthof[c]{{}\leq{}}{{}={}}
			\fillwidthof[c]{\FBE(x_\star)}{\varphi(x_\star)}
		\)
		~iff~
		\(
			x_\star\in\argmin\varphi
		\).
	\end{enumerate}
	In fact, the assumption of twice continuous differentiability of \(f\) can be dropped.
	\begin{proof}
		\begin{proofitemize}
		\item\ref{thm:leq}~
			Since the minimum in \eqref{eq:FBEmin} is attained at \(w=\bar x\), cf. \eqref{eq:Targmin}, we have
			\begin{align*}
			\numberthis
			\label{eq:FBEbarx}
				\FBE(x)
			{}={} &
				f(x)+\innprod{\nabla f(x)}{\bar x-x}+\tfrac{1}{2\gamma}\|\bar x-x\|^2
				{}+{}
				g(\bar x)
			\\
			{}\leq{} &
				f(x)+\innprod{\nabla f(x)}{\bar x-x}+\tfrac{1}{2\gamma}\|\bar x-x\|^2
				{}+{}
				g(x)
				{}+{}
				\innprod{\tfrac1\gamma(x-\bar x)-\nabla f(x)}{\bar x-x}
			\\
			{}={} &
				f(x)
				{}+{}
				g(x)
				{}-{}
				\tfrac{1}{2\gamma}\|x-\bar x\|^2
			\end{align*}
			where in the inequality we used the fact that \(\tfrac1\gamma(x-\bar x)-\nabla f(x)\in\partial g(\bar x)\), cf. \eqref{eq:FBCharacterization}.
		\item\ref{thm:geq}~
			Follows by using \eqref{eq:LipBound} (with \(h=f\), \(u=x\) and \(v=\bar x\)) in \eqref{eq:FBEbarx}.
		\item\ref{thm:=}~
			Follows by \ref{thm:leq} and the optimality condition \eqref{eq:R0}.
		\qedhere
		\end{proofitemize}
	\end{proof}
\end{thm}
Notice that by combining \Cref{thm:geq,thm:leq} we recover the ``sufficient decrease'' condition of (convex) FBS \cite[Thm. 1]{nesterov2013gradient}, that is
\begin{equation}\label{eq:SDFBS}
	\varphi(\bar x)
{}\leq{}
	\varphi(x)
	{}-{}
	\tfrac{2-\gamma L_f}{2\gamma}
	\|x-\bar x\|^2
\end{equation}
holding for all \(x\in\R^n\) with \(\bar x=\T(x)\).

\begin{thm}[FBE: smooth minimization equivalence]\label{thm:equiv}%
	For all \(\gamma>0\)
	\begin{enumerate}
	\item\label{thm:C1}%
		\(\FBE\in C^1(\R^n)\)
		~with~
		\(
			\nabla\FBE
		{}={}
			Q_\gamma\Res
		\).
	\end{enumerate}
	Moreover, if \(\gamma\in(0,\nicefrac{1}{L_f})\) then the following also hold:
	\begin{enumerate}[resume]
	\item\label{thm:infequiv}%
		\(
			\inf\varphi=\inf\FBE
		\);
	\item\label{thm:minimequiv}%
		\(
			x_\star\in\argmin\varphi
		\)
		~iff~
		\(
			x_\star\in\argmin\FBE
		\)
		~iff~
		\(
			\nabla\FBE(x_\star)=0
		\).
	\end{enumerate}
	\begin{proof}
		\begin{proofitemize}
		\item\ref{thm:C1}.
			Since \(f\in C^2(\R^n)\) and \(g^\gamma\in C^1(\R^n)\) (cf. \cref{thm:ME:C1}), from the definition \eqref{eq:FBE} it is apparent that \(\FBE\) is continuously differentiable for all \(\gamma>0\).
			The formula for the gradient was already shown in \eqref{eq:gradFBE}.
		\end{proofitemize}
		Suppose now that \(\gamma<\nicefrac{1}{L_f}\).
		\begin{proofitemize}
		\item\ref{thm:infequiv}.
			\(
				\inf\varphi
			{}\leq{}
				\inf_{x\in\R^n}
				\varphi(\T(x))
			{}\overrel*[\leq]{\ref{thm:geq}}{}
				\inf_{x\in\R^n}
				\FBE(x)
			{}={}
				\inf\FBE
			{}\overrel*[\leq]{\ref{thm:leq}}{}
				\inf\varphi
			\).
		\item\ref{thm:minimequiv}.
			We have
			\[
				x_\star\in\argmin\varphi
			~~\overrel[\Leftrightarrow]{\eqref{eq:R0}}~~
				\Res(x_\star)=0
			~~\Leftrightarrow~~
				Q_\gamma(x_\star)\Res(x_\star)=0
			~~\overrel[\Leftrightarrow]{\ref{thm:C1}}~~
				\nabla\FBE(x_\star)=0,
			\]
			where the second equivalence follows from the invertibility of \(Q_\gamma\).
			
			Suppose now that \(x_\star\in\argmin\FBE\).
			Since \(\FBE\in C^1(\R^n)\) the first-order necessary condition reads \(\nabla\FBE=0\), and from the equivalence proven above we conclude that
			\(
				\argmin\FBE\subseteq\argmin\varphi
			\).
			Conversely, if \(x_\star\in\argmin\varphi\) then
			\[
				\FBE(x_\star)
			{}\overrel*{\ref{thm:=}}{}
				\varphi(x_\star)
			{}={}
				\inf\varphi
			{}\overrel*{\ref{thm:infequiv}}{}
				\inf\FBE,
			\]
			proving \(x_\star\in\argmin\FBE\), hence the inclusion
			\(
				\argmin\FBE\supseteq\argmin\varphi
			\).
		\qedhere
		\end{proofitemize}
	\end{proof}
\end{thm}

\begin{figure}
	\begin{minipage}[t][][t]{0.45\linewidth}
		\vspace*{-\abovecaptionskip}%
		\caption[Forward-backward envelope]{%
			\emph{%
				FBE of the function \(\varphi\) as in \cref{fig:ME} with same parameter \(\gamma=0.2\), relative to the decomposition as the sum of
				\(
					f(x) = x^2+x-1
				\)
				and
				\(
					g(x) = \tfrac13x^3+\indicator_{[0,\infty)}(x)
				\).
			}\newline
			For \(\gamma<\nicefrac{1}{L_f}\) (\(L_f=2\) in this example) at each point \(x\) the FBE \(\FBE\) is the minimum of the quadratic majorization model \(\mathcal M_\gamma^{f,g}({}\cdot{},x)\) for \(\varphi\), the unique minimizer being the proximal gradient point \(\bar x=T_\gamma(x)\).
			The FBE is a differentiable lower bound to \(\varphi\) and since \(f\) is quadratic in this example, it is also smooth and convex (cf. \cref{thm:FBEconvex}).
			In any case, its stationary points and minimizers coincide, and are equivalent to the minimizers of \(\varphi\).
		}%
		\label{fig:FBE}%
	\end{minipage}
	\hfill
	\begin{minipage}[t][][t]{0.5\linewidth}
		\vspace*{0pt}%
		{{%
			\pgfkeys{/pgf/images/include external/.code={\includegraphics[width=\linewidth]{#width=\linewidth}}}%
			\tikzsetnextfilename{FB}%
			\begin{tikzpicture}
\def\gam{0.2}
	\def\xMin{-1.1}
	\def\xMax{2.5}
	\def\yMin{-0.5}
	\def\yMax{6}
	\def\xo{1.8}
	\def\a{1}
	\def\b{-1}
	\def\c{1}
	\def\xstar{max(-\a+sqrt(\a^2-\b), 0)}
	\tikzset{%
		declare function = {%
			f(\x)       = \a*\x^2 + \b*\x + \c;
			gradf(\x)   = 2*\a*\x + \b;
			g(\x)       = (\x < 0) * \YINFTY  +  (\x >= 0) * min(\x^3/3, \YINFTY);
			proxg(\x)   = max(sqrt(1 + 4*\gam*\x) - 1, 0) / (2*\gam); 
			proxphiunconstrained(\x) = (-(2*\gam*\a + 1) + sqrt( (2*\gam*\a + 1)^2 - 4*(\b*\gam-\x)*\gam ) ) / (2*\gam);
			proxphi(\x) = (proxphiunconstrained(\x) >= 0) * proxphiunconstrained(\x);
		},
	}
	\tikzset{%
		declare function = {
			phi(\x)   = (g(\x) >= \YINFTY) * \YINFTY  +  (g(\x) < \YINFTY) * min(f(\x) + g(\x), \YINFTY);
			phiME(\x) = phi(proxphi(\x)) + (\x-proxphi(\x))^2 / (2*\gam);
			gME(\x)   = g(proxg(\x)) + (\x-proxg(\x))^2 / (2*\gam);
			FB(\x)   = proxg(\x-\gam*gradf(\x));
			FBE(\x)  = f(\x) - (\gam/2)*(gradf(\x))^2 + gME(\x-\gam*gradf(\x));
			Q(\x,\w) = (g(\x) >= \YINFTY) * \YINFTY  +  (g(\x) < \YINFTY) * (g(\x) + f(\w) + gradf(\w)*(\x-\w) + (\x-\w)^2/(2*\gam));
			Qphi(\x,\w) = (phi(\x) + (\x-\w)^2/(2*\gam));
		},
	}
	\begin{axis}[myaxis]
		%
		\plotExtFunc[phi]{phi}{\xMin:0.01}
		\addlegendentry{\(\varphi=f+g\)}
		\node[label, anchor=north] at (axis cs: {(\xMin)/2},{\YINFTY}) {\(+\infty\)};
		%
		\addplot[Q, forget plot] {Q(x,\xo)};
		\node[label, color=Q, anchor=north east, outer sep=0pt, font=\tiny] at (axis cs: {0.4*\xo},{Q(0.4*\xo,\xo)}) {\(\mathcal M_\gamma^{f,g}({}\cdot{};x)\)};
		%
		\addplot[Q, dashed, forget plot] {Qphi(x,\xo)};
		\node[label, color=Q, anchor=south west, outer sep=0pt, font=\tiny] at (axis cs: {0.3*\xo},{Qphi(0.3*\xo,\xo)}) {\(\mathcal M_\gamma^\varphi({}\cdot{};x)\)};
		%
		\addplot[ME] {phiME(x)};
		\addlegendentry{\(\varphi^\gamma\)}
		%
		\addplot[FBE] {FBE(x)};
		\addlegendentry{\(\FBE\)}
		%
		\node[nosep] (phi0) at (axis cs: {\xo},{phi(\xo)}) {};
		\drawCoord{phi0}
		\drawDot[color=phi]{phi0};
		\xLabel{phi0}{\(x\vphantom{x_\star}\)}
		\yLabel[color=phi]{phi0}{$\varphi(x)$}
		%
		\node[nosep] (me0) at (axis cs: {\xo},{phiME(\xo)}) {};
		\drawCoord{me0}
		\drawDot[color=ME]{me0};
		\yLabel[color=ME]{me0}{\(\varphi^\gamma(x)\)}
		%
		\node[nosep] (fbe0) at (axis cs: {\xo},{FBE(\xo)}) {};
		\drawCoord{fbe0}
		\drawDot[color=FBE]{fbe0}
		\yLabel[color=FBE]{fbe0}{\(\FBE(x)\)}
		%
		\node[nosep] (Q0) at (axis cs: {FB(\xo)},{FBE(\xo)}) {};
		\drawCoord{Q0}
		\drawDot[color=Q]{Q0}
		\xLabel{Q0}{\(\bar x\vphantom{x_\star}\)}
		%
		\node[nosep] (phiP0) at (axis cs: {FB(\xo)},{phi(FB(\xo))}) {};
		\drawCoord{phiP0}
		\drawDot[color=phi]{phiP0}
		\yLabel[color=phi]{phiP0}{\(\varphi(\bar x)\)}
		%
		\node[nosep] (xmin) at (axis cs: {\xstar},{phi(\xstar)}) {};
		\drawCoord[dashed][dashed]{xmin}
		\drawDot{xmin}
		\xLabel[anchor=south]{xmin}{\(x_\star\vphantom{x_\star}\)}
		\yLabel[anchor=north east]{xmin}{\(\min\varphi=\min\FBE\)}
	\end{axis}
\end{tikzpicture}%
		}}%
	\end{minipage}
\end{figure}

\begin{prop}[FBE and Moreau envelope {\cite[Thm. 2]{nesterov2013gradient}}]\label{thm:FBE-Moreau}%
	For any \(\gamma\in(0,\nicefrac{1}{L_f})\), it holds that
	\[
		\varphi^{\frac{\gamma}{1-\gamma L_f}}
	{}\leq{}
		\varphi_\gamma
	{}\leq{}
		\varphi^{\gamma}.
	\]
	\begin{proof}
		We have
		\begin{align*}
			\FBE(x)
		{}={} &
			\min_{w\in\R^n}\set{
				\smashunderbracket{f(x)+\innprod{\nabla f(x)}{w-x}}{}
				{}+{}
				\tfrac{1}{2\gamma}\|w-x\|^2
				{}+{}
				g(w)
			}
		\\
		{}\overrel[\leq]{\eqref{eq:LipBound}}{} &
			\min_{w\in\R^n}\set{
				\smashoverbracket{
					\fillwidthof[c]{f(x)+\innprod{\nabla f(x)}{w-x}}{f(w)-\tfrac{L_f}{2}\|w-x\|^2}
				}{}
				{}+{}
				\tfrac{1}{2\gamma}\|w-x\|^2
				{}+{}
				g(w)
			}
		\\
		{}={} &
			\min_{w\in\R^n}\set{
				f(w)
				{}+{}
				g(w)
				{}+{}
				\tfrac{1-\gamma L_f}{2\gamma}\|w-x\|^2
			}
		{}={}
			\varphi^{\frac{\gamma}{1-\gamma L_f}}(x).
		\end{align*}
		Using the upper bound in \eqref{eq:LipBound} instead, yields the other inequality.
	\end{proof}
\end{prop}
Since \(\FBE\) is upper bounded by the \(\gamma^{-1}\)-smooth function \(\varphi^\gamma\) with which it shares the set of minimizers \(\mathcal X_\star\), from \eqref{eq:LipBound} we easily infer the following quadratic upper bound.
\begin{cor}[Global quadratic upper bound]\label{thm:QUB}%
	If \(\mathcal X_\star\neq\emptyset\), then
	\[
		\FBE(x)-\varphi_\star
	{}\leq{}
		\tfrac{1}{2\gamma}
		\dist(x,\mathcal X_\star)^2
		\quad
		\forall x\in\R^n.
	\]
\end{cor}

Although the FBE may fail to be convex, for \(\gamma<\nicefrac{1}{L_f}\) its stationary points and minimizers coincide and are the same as those of the original function \(\varphi\).
That is, the minimization of \(\varphi\) is equivalent to the minimization of the differentiable function \(\FBE\).
This is a clear analogy with the Moreau envelope, which in fact is the special case of the FBE corresponding to \(f\equiv 0\) in the decomposition of \(\varphi\).
In the next result we tighten the claims of \Cref{thm:C1} when \(f\) is a convex quadratic function, showing that in this case the FBE is convex and smooth and thus recover all the properties of the Moreau envelope.
\begin{thm}[FBE: convexity \& smoothness for quadratic \(f\) {\cite[Prop. 4.4]{giselsson2018envelope}}]\label{thm:FBEconvex}
	Suppose that \(f\) is convex quadratic, namely \(f(x)=\tfrac12\innprod{x}{Hx}+\innprod hx\) for some \(H\in\sym_+(\R^n)\) and \(h\in\R^n\).
	Then, for all \(\gamma\in(0,\nicefrac{1}{L_f}]\) the FBE \(\FBE\) is convex and smooth, with
	\[
		L_{\FBE}
	{}={}
		\tfrac{1-\gamma\mu_f}{\gamma}
	\quad\text{and}\quad
		\mu_{\FBE}
	{}={}
		\min\set{
			\mu_f(1-\gamma\mu_f),
			L_f(1-\gamma L_f)
		},
	\]
	where \(L_f=\lambda_{\rm max}(H)\) and \(\mu_f=\lambda_{\rm min}(H)\).
	In particular, when \(f\) is \(\mu_f\)-strongly convex the strong convexity of \(\FBE\) is maximized for \(\gamma=\frac{1}{\mu_f+L_f}\), in which case
	\[
		L_{\FBE}
	{}={}
		L_f
	\quad\text{and}\quad
		\mu_{\FBE}
	{}={}
		\tfrac{L_f\mu_f}{\mu_f+L_f}.
	\]
	\begin{proof}
		Letting \(Q\coloneqq\I-\gamma H\), we have that \(Q_\gamma\equiv Q\) and \(\Fw x=Qx-\gamma h\).
		Therefore,
		\begin{align*}
			\gamma
			\innprod{
				\nabla\FBE(x)-\nabla\FBE(y)
			}{
				x-y
			}
		{}\overrel*{\eqref{eq:gradFBE}}{} &
			\innprod{
				Q(\Res(x)-\Res(y))
			}{
				x-y
			}
		\\
		{}={} &
			\innprod{
				Q(x-y)
			}{
				x-y
			}
			{}-{}
			\innprod{
				Q(\T(x)-\T(y))
			}{
				x-y
			}
		\\
		{}={} &
			\|x-y\|_Q^2
		\\
		&
			{}-{}
			\innprod{
				\prox_{\gamma g}(Qx-\gamma h)
				{}-{}
				\prox_{\gamma g}(Qy-\gamma h)
			}{
				Q(x-y)
			}.
		\end{align*}
		From the \emph{firm} nonexpansiveness of \(\prox_{\gamma g}\) (see \cite[Prop.s 4.35(iii) and 12.28]{bauschke2017convex}) it follows that
		\[
			0
		{}\leq{}
			\innprod{
				\prox_{\gamma g}(Qx-\gamma h)
				{}-{}
				\prox_{\gamma g}(Qy-\gamma h)
			}{
				Q(x-y)
			}
		{}\leq{}
			\|Q(x-y)\|^2.
		\]
		By combining with the previous inequality, we obtain
		\[
			\tfrac1\gamma\|x-y\|_{Q-Q^2}^2
		{}\leq{}
			\innprod{
				\nabla\FBE(x)-\nabla\FBE(y)
			}{
				x-y
			}
		{}\leq{}
			\tfrac1\gamma
			\|x-y\|_Q^2.
		\]
		Since \(\lambda_{\rm min}(Q)=1-\gamma L_f\) and \(\lambda_{\rm max}(Q)=1-\gamma\mu_f\), from \cref{thm:eigs} we conclude that
		\[
			\mu_{\FBE}
			\|x-y\|^2
		{}\leq{}
			\innprod{
				\nabla\FBE(x)-\nabla\FBE(y)
			}{
				x-y
			}
		{}\leq{}
			L_{\FBE}
			\|x-y\|^2
		\]
		with \(\mu_{\FBE}\) and \(L_{\FBE}\) as in the statement, hence the claim, cf. \cref{eq:Lmu}.
%
%
	\end{proof}
\end{thm}

\begin{lem}\label{thm:FBE_R2}%
	Suppose that \(\varphi\) has the quadratic growth with constants \((\mu,\nu)\), and let \(\varphi_\star\coloneqq\min\varphi\).
	Then, for all \(\gamma\in(0,\nicefrac{1}{L_f}]\) and \(x\in\lev_{\leq\varphi_\star+\nu}\FBE\) it holds that
	
	\[
		\varphi_\gamma(x)-\varphi_\star
	{}\leq{}
		\gamma
		\bigl[
			\tfrac12
			{}+{}
			(1+\nicefrac{2}{\gamma\mu})
			(1+\gamma L_f)^2
		\bigr]
		\|\Res(x)\|^2.
	\]
	\begin{proof}
		Fix \(x\in\lev_{\leq\varphi_\star+\nu}\FBE\) and let \(\bar x\coloneqq\T(x)\).
		We have
		\begin{align*}
			\FBE(x)-\varphi_\star
		{}\overrel[\leq]{\ref{thm:geq}}{} &
			\tfrac\gamma2
			\|\Res(x)\|^2
			{}+{}
			\varphi(\bar x)-\varphi_\star
		\\
		{}\overrel[\leq]{\ref{thm:cost_dd}}{} &
			\tfrac\gamma2
			\|\Res(x)\|^2
			{}+{}
			\dist(\bar x,\mathcal X_\star)
			\dist(0,\partial\varphi(\bar x))
		\\
		{}\overrel[\leq]{\ref{thm:r_d}}{} &
			\bigl[
				\tfrac\gamma2
				\|\Res(x)\|
				{}+{}
				(1+\gamma L_f)
				\dist(\bar x,\mathcal X_\star)
			\bigr]
			\|\Res(x)\|
		\shortintertext{%
			and since \(\bar x\in\lev_{\leq\varphi_\star+\nu}\varphi\) (cf. \cref{thm:geq}), from \cref{thm:EB:Res} we can bound the quantity \(\dist(\bar x,\mathcal X_\star)\) in terms of the residual as
		}
		{}\leq{} &
			\bigl[
				\tfrac\gamma2
				\|\Res(x)\|
				{}+{}
				(\gamma+\nicefrac2\mu)(1+\gamma L_f)^2
				\|\Res(\bar x)\|
			\bigr]
			\|\Res(x)\|.
		\end{align*}
		The proof now follows from the inequality \(\|\Res(\bar x)\|\leq\|\Res(x)\|\), see \cite[Thm. 10.12]{beck2017first}, after easy algebraic manipulations.
	\end{proof}
\end{lem}

\subsection{Further equivalence properties}
\begin{prop}[Equivalence of level boundedness]\label{thm:LB}%
	For any \(\gamma\in(0,\nicefrac{1}{L_f})\), \(\varphi\) has bounded level sets iff \(\FBE\) does.
	\begin{proof}
		\cref{thm:sandwich} implies that
		\(
			\lev_{\leq\alpha}\varphi
		\subseteq
			\lev_{\leq\alpha}\FBE
		\)
		for all \(\alpha\in\R\), therefore level boundedness of \(\FBE\) implies that of \(\varphi\).
		Conversely, suppose that \(\FBE\) is not level bounded, and consider \(\seq{x_k}\subseteq\lev_{\leq\alpha}\FBE\) with \(\|x_k\|\to\infty\).
		Then from \cref{thm:sandwich} it follows that
		\(
			\varphi(\bar x_k)
		{}\leq{}
			\FBE(x_k)-\tfrac{1}{2\gamma}\|x_k-\bar x_k\|^2
		{}\leq{}
			\alpha-\tfrac{1}{2\gamma}\|x_k-\bar x_k\|^2
		\),
		where \(\bar x_k=\T(x_k)\).
		In particular, \(\seq{\bar x_k}\subseteq\lev_{\leq\alpha}\varphi\).
		If \(\seq{\bar x_k}\) is bounded, then \(\inf\varphi=-\infty\); otherwise, \(\lev_{\leq\alpha}\varphi\) contains the unbounded sequence \(\seq{\bar x_k}\).
		Either way, \(\varphi\) cannot be level bounded.
	\end{proof}
\end{prop}

\begin{prop}[Equivalence of quadratic growth]\label{thm:QGequiv}%
	Let \(\gamma\in(0,\nicefrac{1}{L_f})\) be fixed.
	Then,
	\begin{enumerate}
	\item\label{thm:QGequiv:phi}
		if \(\varphi\) satisfies the quadratic growth condition with constants \((\mu,\nu)\), then so does \(\FBE\) with constants \((\mu',\nu)\), where
		\(
			\mu'
		{}\coloneqq{}
			\tfrac{1-\gamma L_f}{(1+\gamma L_f)^2}
			\tfrac{\mu\gamma}{(2+\gamma\mu)^2}
			\mu
		\);
	\item\label{thm:QGequiv:FBE}
		conversely, if \(\FBE\) satisfies the quadratic growth condition, then so does \(\varphi\) with same constants.
	\end{enumerate}
	\begin{proof}
		Since \(\varphi\) and \(\FBE\) have same infimum and minimizers (cf. \cref{thm:equiv}), \ref{thm:QGequiv:FBE} is a straightforward consequence of the fact that \(\FBE\leq\varphi\) (cf. \cref{thm:leq}).
		
		Conversely, suppose that \(\varphi\) satisfies the quadratic growth with constants \((\mu,\nu)\).
		Then, for all
		\(
			x
		{}\in{}
			\lev_{\leq\varphi_\star+\nu}\FBE
		\)
		we have that
		\(
			\bar x\coloneqq\T(x)
		{}\in{}
			\lev_{\leq\varphi_\star+\nu}\varphi
		\),
		therefore
		\[
			\FBE(x)-\varphi_\star
		{}\overrel[\geq]{\ref{thm:geq}}{}
			\varphi(\bar x)-\varphi_\star
			{}+{}
			\gamma\tfrac{1-\gamma L_f}{2}\|\Res(x)\|^2
		{}\geq{}
			\tfrac{\mu'}{2}
			\dist(x,\mathcal X_\star),
		\]
		where in the last inequality we discarded the term \(\varphi(\bar x)-\varphi_\star\geq0\) and used \cref{thm:EB:Res} to lower bound \(\|\Res(x)\|^2\).
	\end{proof}
\end{prop}

\begin{cor}[Equivalence of strong minimality]\label{thm:equiv:strmin}%
	For all \(\gamma\in(0,\nicefrac{1}{L_f})\), a point \(x_\star\) is a (locally) strong minimizer for \(\varphi\) iff it is a (locally) strong minimizer for \(\FBE\).
\end{cor}

Lastly, having showed that for convex functions the quadratic growth can be extended to arbitrary level sets (cf. \cref{thm:QGglobal}), an interesting consequence of \Cref{thm:QGequiv} is that, although \(\FBE\) may fail to be convex, it enjoys the same property.
\begin{cor}[FBE: globality of quadratic growth]%
	Let \(\gamma\in(0,\nicefrac{1}{L_f})\) and suppose that \(\FBE\) satisfies the quadratic growth with constants \((\mu,\nu)\).
	Then, for every \(\nu'>\nu\) there exists \(\mu'>0\) such that \(\FBE\) satisfies the quadratic growth with constants \((\mu',\nu')\).
\end{cor}

		\subsection{Second-order properties}\label{sec:2nd}%
				Although \(\FBE\) is continuously differentiable over \(\R^n\), it fails to be \(C^2\) in most cases; since \(g\) is nonsmooth, its Moreau envelope \(g^\gamma\) is hardly ever \(C^2\).
	For example, if \(g\) is real valued then \(g^\gamma\) is  \(C^2\) (and \(\prox_{\gamma g}\) is \(C^1\)) if and only if \(g\) is \(C^2\) \cite{lemarechal1997practical}.
	Therefore, we hardly ever have the luxury of assuming continuous differentiability of \(\nabla \FBE\) and we must resort to generalized notions of differentiability stemming from nonsmooth analysis.
	Specifically, our analysis is largely based on generalized differentiability properties of \(\prox_{\gamma g}\) which we study next.
	\begin{thm}\label{th:JacProx}%
		For all \(x\in\R^n\), \(\partial_C(\prox_{\gamma g})(x)\neq\emptyset\) and any \(P\in\partial_C(\prox_{\gamma g})(x)\) is a symmetric positive semidefinite matrix that satisfies \(\|P\|\leq 1\).
		\begin{proof}
			Nonempty-valuedness of \(\partial_C(\prox_{\gamma g})\) is due to Lipschitz continuity of \(\prox_{\gamma g}\).
			Moreover, since \(g\) is convex, its Moreau envelope is a convex function as well, therefore every element of \(\partial_C(\nabla g^\gamma)(x)\) is a symmetric positive semidefinite matrix (see \eg \cite[\S8.3.3]{facchinei2003finite}).
			Due to \cref{thm:ME:C1}, we have that
			\(
				\prox_{\gamma g}(x)
			{}={}
				x-\gamma\nabla g^\gamma(x)
			\),
			therefore
			\begin{equation}\label{eq:JacProx}
				\partial_C(\prox_{\gamma g})(x)
			{}={}
				\I-\gamma\partial_C(\nabla g^\gamma)(x).
			\end{equation}
			The last relation holds with equality (as opposed to inclusion in the general case) due to the fact that one of the summands is continuously differentiable.
			Now, from \eqref{eq:JacProx} we easily infer that every element of \(\partial_C(\prox_{\gamma g})(x)\) is a symmetric matrix.
			Since \(\nabla g^\gamma(x)\) is Lipschitz continuous with Lipschitz constant \(\gamma^{-1}\), using \cite[Prop. 2.6.2(d)]{clarke1990optimization}, we infer that every \(H\in\partial_C(\nabla g^\gamma)(x)\) satisfies \(\|H\|\leq\gamma^{-1}\).
			Now, according to \eqref{eq:JacProx} it holds that
			\[
				P\in\partial_C(\prox_{\gamma g})(x)
			\quad\Leftrightarrow\quad
				P=I-\gamma H,
			\quad
				H\in\partial_C(\nabla g^\gamma)(x).
			\]
			Therefore, for every \(d\in\R^n\) and \(P\in\partial_C(\prox_{\gamma g})(x)\),
			\[
				\innprod{d}{Pd}
			{}={}
				\|d\|^2-\gamma\innprod{d}{Hd}
			{}\geq{}
				\|d\|^2-\gamma\gamma^{-1}\|d\|^2
			{}={}
				0.
			\]
			On the other hand, since \(\prox_{\gamma g}\) is Lipschitz continuous with Lipschitz constant 1, using \cite[Prop. 2.6.2(d)]{clarke1990optimization} we obtain that \(\|P\|\leq 1\) for all \(P\in\partial_C(\prox_{\gamma g})(x)\).
		\end{proof}
	\end{thm}

	We are now in a position to construct a generalized Hessian for \(\FBE\) that will allow the development of Newton-like methods with fast asymptotic convergence rates.
	An obvious route to follow would be to assume that \(\nabla\FBE\) is semismooth and employ \(\partial_C(\nabla\FBE)\) as a generalized Hessian for \(\FBE\).
	However, this approach would require extra assumptions on \(f\) and involve complicated operations to evaluate elements of \(\partial_C(\nabla\FBE)\).
	On the other hand, what is really needed to devise Newton-like algorithms with fast local convergence rates is a \emph{linear Newton approximation (LNA)}, cf. \Cref{def:LNA}, at some stationary point of \(\FBE\), which by \Cref{thm:minimequiv} is also a minimizer of \(\varphi\), provided that \(\gamma\in(0,\nicefrac{1}{L_f})\).

	The approach we follow is largely based on \cite{sun1997computable}, \cite[Prop. 10.4.4]{facchinei2003finite}.
	Without any additional assumptions we can define a set-valued mapping \(\ffunc{\hat\partial^2\varphi_\gamma}{\R^n}{\R^{n\times n}}\) with full domain and whose elements have a simpler form than those of \(\partial_C(\nabla\FBE)\), which serves as a LNA for \(\nabla\FBE\) at any stationary point \(x_\star\) provided \(\prox_{\gamma g}\) is semismooth at \(\Fw{x_\star}\).
	We call it \DEF{approximate generalized Hessian} of \(\FBE\) and it is given by
	\begin{equation}\label{eq:GenHess}
		\hat\partial^2\FBE(x)
	\coloneqq
		\set{
			\gamma^{-1}Q_\gamma(x)(\I-PQ_\gamma(x))
		}[
			P\in\partial_C(\prox_{\gamma g})(\Fw x)
		].
	\end{equation}
	Notice that if \(f\) is quadratic, then \(\hat\partial^2\FBE\equiv\partial_C\nabla\FBE\); more generally, the key idea in the definition of \(\hat\partial^2\FBE\), reminiscent of the Gauss-Newton method for nonlinear least-squares problems, is to omit terms vanishing at \(x_\star\) that contain third-order derivatives of \(f\).	
	
	\begin{prop}\label{thm:LNAR}%
		Let \(\bar x\in\R^n\) and \(\gamma>0\) be fixed.
		If \(\prox_{\gamma g}\) is (\(\vartheta\)-order) semismooth at \(\Fw{\bar x}\) (and \(\nabla^2f\) is \(\vartheta\)-H\"older continuous around \(\bar x\)), then
		\begin{equation}\label{eq:JR}
			\mathcal R_\gamma(x)
		{}\coloneqq{}
			\set{\gamma^{-1}(\I-PQ_\gamma(x))}[
				P\in\partial_C\FB x
			]
		\end{equation}
		is a (\(\vartheta\)-order) LNA for \(\Res\) at \(\bar x\).
		\begin{proof}
			We shall prove only the \(\vartheta\)-order semismooth case, as the other one is shown by simply replacing all occurrences of \(O(\|{}\cdot{}\|^{1+\vartheta})\) with \(o(\|{}\cdot{}\|)\) in the proof.
			Let \(q_\gamma=\Fw{}\) be the forward operator, so that the forward-backward operator \(\T\) can be expressed as \(\T=\prox_{\gamma g}\circ q_\gamma\).
			With a straightforward adaptation of the proof of \cite[Prop. 7.2.9]{facchinei2003finite} to include the \(\vartheta\)-H\"olderian case, it can be shown that
			\begin{equation}\label{eq:FPtheta}
				q_\gamma(x)-q_\gamma(\bar x)-Q_\gamma(x)(x-\bar x)
			{}={}
				O(\|x-\bar x\|^{1+\vartheta}).
			\end{equation}
			Moreover, since \(\nabla f\) is Lipschitz continuous and thus so is \(q_\gamma\), we also have
			\begin{equation}\label{eq:qLip}
				q_\gamma(x)-q_\gamma(\bar x)
			{}={}
				O(\|x-\bar x\|).
			\end{equation}
			Let \(U_x\in\mathcal\Res(x)\) be arbitrary; then, there exists \(P_x\in\partial_C\FB x\) such that
			\(
				U_x
			{}={}
				\gamma^{-1}(\I-P_xQ_\gamma(x))(\bar x-x)
			\).
			We have
			\begin{align*}
			&
				\Res(x)
				{}+{}
				U_x(\bar x-x)
				{}-{}
				\Res(\bar x)
			\\
			{}={} &
				\Res(x)
				{}+{}
				\gamma^{-1}(\I-P_xQ_\gamma(x))(\bar x-x)
				{}-{}
				\Res(\bar x)
			\\
			{}={} &
				\gamma^{-1}
				\bigl[
					\prox_{\gamma g}(q_\gamma(\bar x))
					{}-{}
					\prox_{\gamma g}(q_\gamma(x))
					{}-{}
					P_xQ_\gamma(x)(\bar x-x)
				\bigr]
			\shortintertext{%
				due to \(\vartheta\)-order semismoothness of \(\prox_{\gamma g}\) at \(q_\gamma(\bar x)\),
			}
			{}={} &
				\gamma^{-1}
				P_x
				\bigl[
					q_\gamma(\bar x)-q_\gamma(x)
					{}+{}
					O(\|q_\gamma(\bar x)-q_\gamma(x)\|^{1+\vartheta})
					{}-{}
					Q_\gamma(x)
					(\bar x-x)
				\bigr]
			\\
			{}\overrel{\eqref{eq:qLip}}{} &
				\gamma^{-1}
				P_x
				\bigl[
					q_\gamma(\bar x)-q_\gamma(x)
					{}-{}
					Q_\gamma(\bar x)
					(\bar x-x)
					{}+{}
					O(\|\bar x-x\|^{1+\vartheta})
				\bigr]
			\\
			{}\overrel{\eqref{eq:FPtheta}}{} &
				\gamma^{-1}
				P_x
				O(\|\bar x-x\|^{1+\vartheta})
			{}={}
				O(\|\bar x-x\|^{1+\vartheta}),
			\end{align*}
			where in the last equality we used the fact that \(\|P_x\|\leq1\), cf. \cref{th:JacProx}.
%
%
%
		\end{proof}
	\end{prop}
	
	\begin{cor}\label{prop:LNAprops1}%
		Let \(\gamma\in(0,1/L_f)\) and \(x_\star\in\mathcal X_\star\).
		If \(\prox_{\gamma g}\) is (\(\vartheta\)-order) semismooth at \(\Fw{x_\star}\) (and \(\nabla^2f\) is locally \(\vartheta\)-H\"older continuous around \(x_\star\)), then \(\hat\partial^2\FBE\) is a (\(\vartheta\)-order) LNA for \(\nabla\FBE\) at \(x_\star\).
		\begin{proof}
			Let
			\(
				H_x
			{}\in{}
				\hat\partial^2\FBE(x)
			{}={}
				\set{Q_\gamma(x)U}[U\in\mathcal R_\gamma(x)]
			\),
			so that \(H_x=Q_\gamma(x)U_x\) for some \(U_x\in\mathcal R_\gamma(x)\).
			Then,
			\begin{align*}
				\|
					\nabla\FBE(x)+H_x(x_\star-x)-\nabla\FBE(x_\star)
				\|
			{}={} &
				\|
					Q_\gamma(x)R_\gamma(x)+Q_\gamma(x)U_x(x-x_\star))
				\|
			\\
			{}={} &
				\|
					Q_\gamma(x)[R_\gamma(x)+U_x(x-x_\star)-R_\gamma(x_\star)]
				\|
			\\
			{}\leq{} &
				\|R_\gamma(x)+U_x(x-x_\star)-R_\gamma(x_\star)\|,
			\end{align*}
			where in the equalities we used the fact that \(\nabla\FBE(x_\star)=R_\gamma(x_\star)=0\), and in the inequality the fact that \(\|Q_\gamma\|\leq1\).
			Since \(\mathcal R_\gamma\) is a (\(\vartheta\)-order) LNA of \(R_\gamma\) at \(x_\star\), the last term is \(o(\|x-x_\star\|)\) (resp. \(O(\|x-x_\star\|^{1+\vartheta})\)).
		\end{proof}
	\end{cor}
	As shown in the next result, although the FBE is in general not convex, for \(\gamma\) small enough every element of \(\hat\partial^2\FBE(x)\) is a (symmetric and) positive semidefinite matrix.
	Moreover, the eigenvalues are lower and upper bounded uniformly over all \(x\in\R^n\).
	\begin{prop}\label{prop:PSDHess}
		Let \(\gamma\leq\nicefrac{1}{L_f}\) and \(H\in\hat\partial^2\FBE(x)\) be fixed.
		Then, \(H\in\sym_+(\R^n)\) with
		\[
			\lambda_{\rm min}(H)
			{}={}
			\min\set{(1-\gamma\mu_f)\mu_f,(1-\gamma L_f)L_f}
		~~\text{and}~~
			\lambda_{\rm max}(H)
			{}={}
			\gamma^{-1}(1-\gamma\mu_f),
		\]
		where \(\mu_f\geq 0\) is the modulus of strong convexity of \(f\).
		\begin{proof}
			Fix \(x\in\R^n\) and let \(Q\coloneqq Q_\gamma(x)\).
			Any \(H\in\hat\partial^2\FBE(x)\) can be expressed as
			\(
				H
			{}={}
				\gamma^{-1}
				Q(\I-PQ)
			\)
			for some \(P\in\partial_C(\prox_{\gamma g})(\Fw x)\).
			Since both \(Q\) and \(P\) are symmetric (cf. \cref{th:JacProx}), it follows that so is \(H\).
			Moreover, for all \(d\in\R^n\)
			\begin{qedalign*}
			\numberthis\label{eq:QPQ}
				\innprod{Hd}{d}
			{}={} &
				\gamma^{-1}\innprod{Qd}{d}
				{}-{}
				\gamma^{-1}\innprod{PQd}{Qd}
			\\
			{}\overrel[\geq]{\ref{th:JacProx}}{} &
				\gamma^{-1}\innprod{Qd}{d}
				{}-{}
				\gamma^{-1}\|Qd\|^2
			\\
			{}={} &
				\innprod{(\I-\gamma\nabla^2f(x))\nabla^2f(x)d}{d}
			\\
			{}\overrel[\geq]{\ref{thm:eigs}}{} &
				\min\set{
					(1-\gamma\mu_f)\mu_f,\,
					(1-\gamma L_f)L_f
				}
				\|d\|^2.
			\shortintertext{%
				On the other hand, since \(P\succeq0\) (cf. \cref{th:JacProx}) and thus \(QPQ\succeq0\), we can upper bound \eqref{eq:QPQ} as
			}
				\innprod{Hd}{d}
			{}\leq{} &
				\gamma^{-1}\innprod{Qd}{d}
			{}\leq{}
				\|Q\|\|d\|^2
			{}\leq{}
				\gamma^{-1}(1-\gamma\mu_f)\|d\|^2.
			\end{qedalign*}%
		\end{proof}
	\end{prop}
	The next lemma links the behavior of the FBE close to a solution of \eqref{eq:GenProb} and a nonsingularity assumption on the elements of \(\hat\partial^2\FBE(x_\star)\).
	Part of the statement is similar to \cite[Lem. 7.2.10]{facchinei2003finite}; however, here \(\nabla\FBE\) is not required to be locally Lipschitz around \(x_\star\).
	\begin{lem}\label{lem:sharpMin}
		Let \(x_\star\in\argmin\varphi\) and \(\gamma\in(0,1/L_f)\).
		If \(\prox_{\gamma g}\) is semismooth at \(\Fw{x_\star}\), then the following conditions are equivalent:
		\begin{enumerateq}
		\item\label{thm:sharpMin:strmin}%
			\(x_\star\) is a locally strong minimum for \(\varphi\) (or, equivalently, for \(\FBE\));
		\item\label{thm:sharpMin:nonsing}%
			every element of \(\hat\partial^2\FBE(x_\star)\) is nonsingular.
		\end{enumerateq}
		In any such case, there exist \(\delta,\kappa>0\) such that 
		\[
			\|x-x_\star\|
		{}\leq{}
			\kappa\|R_\gamma(x)\|
		~~\text{and}~~
			\max\set{\|H\|,\|H^{-1}\|}
		{}\leq{}
			\kappa,
		\]
		for any
		\(
			x\in\ball{x_\star}{\delta}
		\)
		and
		\(
			H\in\hat\partial^2\FBE(x)
		\).
		\begin{proof}
			Observe first that \cref{prop:LNAprops1} ensures that \(\hat\partial^2\FBE\) is a LNA of \(\nabla\FBE\) at \(x_\star\), thus semicontinuous and compact valued (by definition).
			In particular, the last claim follows from \cite[Lem. 7.5.2]{facchinei2003finite}.
			\begin{proofitemize}
			\item\ref{thm:sharpMin:strmin}~\(\Rightarrow\)~\ref{thm:sharpMin:nonsing}~~%
				It follows from \cref{thm:equiv:strmin} that there exists \(\mu,\delta>0\) such that
				\(
					\FBE(x)-\varphi_\star
				{}\geq{}
					\tfrac\mu2\|x-x_\star\|^2
				\)
				for all \(x\in\ball{x_\star}{\delta}\).
				In particular, for all \(H\in\hat\partial^2\FBE(x_\star)\) and \(x\in\ball{x_\star}{\delta}\) we have
				\[
					\tfrac\mu2\|x-x_\star\|^2
				{}\leq{}
					\FBE(x)-\varphi_\star
				{}={}
					\tfrac12\innprod{H(x-x_\star)}{x-x_\star}
					{}+{}
					o(\|x-x_\star\|^2).
				\]
				Let \(v_{\rm min}\) be a unitary eigenvector of \(H\) corresponding to the minimum eigenvalue \(\lambda_{\rm min}(H)\).
				Then, for all \(\varepsilon\in(-\delta,\delta)\) the point \(x_\varepsilon=x_\star+\varepsilon v_{\rm min}\) is \(\delta\)-close to \(x_\star\) and thus
				\[
					\tfrac12\lambda_{\rm min}(H)\varepsilon^2
				{}\geq{}
					\tfrac\mu2\varepsilon^2
					{}+{}
					o(\varepsilon^2)
				{}\geq{}
					\tfrac\mu4\varepsilon^2,
				\]
				where the last inequality holds up to possibly restricting \(\delta\) (and thus \(\varepsilon\)).
				The claim now follows from the arbitrarity of \(H\in\hat\partial^2\FBE(x_\star)\).
			\item\ref{thm:sharpMin:strmin}~\(\Leftarrow\)~\ref{thm:sharpMin:nonsing}~~%
				Easily follows by reversing the arguments of the other implication.
				\qedhere
			\end{proofitemize}
		\end{proof}
	\end{lem}

	\section{Forward-backward truncated-Newton algorithm (FBTN)}\label{sec:Algorithm}
		Having established the equivalence between minimizing \(\varphi\) and \(\FBE\), we may recast problem \eqref{eq:GenProb} into the smooth unconstrained minimization of the FBE.
Under some assumptions the elements of \(\hat\partial^2\FBE\) mimick second-order derivatives of \(\FBE\), suggesting the employment of Newton-like update directions \(d=-(H+\delta\I)^{-1}\nabla\FBE(x)\) with \(H\in\hat\partial^2\FBE(x)\) and \(\delta>0\) (the regularization term \(\delta\I\) ensures the well definedness of \(d\), as \(H\) is positive semidefinite, see \cref{prop:PSDHess}).
If \(\delta\) and \(\varepsilon\) are suitably selected, under some nondegeneracy assumptions updates \(x^+=x+d\) are locally superlinearly convergent.
Since such \(d\)'s are directions of descent for \(\FBE\), a possible globalization strategy is an Armijo-type linesearch.
Here, however, we follow the simpler approach proposed in \cite{stella2017simple,themelis2018forward} that exploits the basic properties of the FBE investigated in \Cref{sec:Basic}.
As we will discuss shortly after, this is also advantageous from a computational point of view, as it allows an arbitrary warm starting for solving the underlying linear system.

Let us elaborate on the linesearch.
To this end, let \(x\) be the current iterate; then, \cref{thm:sandwich} ensures that \(\FBE(\T(x))\leq\FBE(x)-\gamma\frac{1-\gamma L_f}{2}\|\Res(x)\|^2\).
Therefore, unless \(\Res(x)=0\), in which case \(x\) would be a solution, for any \(\sigma\in(0,\gamma\frac{1-\gamma L_f}{2})\) the strict inequality \(\FBE(\T(x))<\FBE(x)-\sigma\|\Res(x)\|^2\) is satisfied.
Due to the continuity of \(\FBE\), all points sufficiently close to \(\T(x)\) will also satisfy the inequality, thus so will the point \(x^+=(1-\tau)\T(x)+\tau(x+d)\) for small enough stepsizes \(\tau\).
This fact can be use to enforce the iterates to \emph{sufficiently} decrease the value of the FBE, cf. \eqref{eq:LS}, which straightforwardly implies optimality of all accumulation points of the generated sequence.
We defer the details to the proof of \Cref{thm:subseq}.
In \Cref{thm:tau1,thm:superlinear} we will provide conditions ensuring acceptance of unit stepsizes so that the scheme reduces to a regularized version of the (undamped) linear Newton method \cite[Alg. 7.5.14]{facchinei2003finite} for solving \(\nabla\FBE(x)=0\), which, under due assumptions, converges superlinearly.
\begin{algorithm}[t]
	\algcaption{%
		({\sf FBTN}) Forward-Backward Truncated-Newton method%
	}%
	\label{alg:FBTN}%
	\begin{algorithmic}[1]
	\Require
		\begin{tabular}[t]{@{}l@{}}
			\(\gamma\in(0,\nicefrac{1}{L_f})\);~
			\(\sigma\in(0,\frac{\gamma(1-\gamma L_f)}{2})\);~
			\(\bar\eta,\zeta\in(0,1)\);~
			\(\rho,\nu\in(0,1]\)
		\\
			initial point \(x_0\in\R^n\);~
			accuracy \(\varepsilon>0\)
		\end{tabular}
	\Provide
		\(\varepsilon\)-suboptimal solution \(x^k\) (\ie such that \(\|\Res(x^k)\|\leq\varepsilon\))
	\Initialize
		\(k\gets 0\)
	\While{~\(\|\Res(x^k)\|>\varepsilon\)~}\label{step:FBN:1}%
		\State
			\(
				\delta_k\gets\zeta\|\nabla\varphi_\gamma(x^k)\|^\nu
			\),~
			\(
				\eta_k\gets\min\set{\bar\eta,\,\|\nabla\varphi_\gamma(x^k)\|^\rho}
			\),~
			\(
				\varepsilon_k\gets\eta_k\|\nabla\varphi_\gamma(x^k)\|
			\)
		\State\label{step:FBN:CG}%
			Apply \refCG\ to find an \(\varepsilon_k\)-approximate solution \(d^k\) to
			\[
				\bigl[H_k+\delta_k\I\bigr]d^k
			{}\approx{}
				-\nabla\varphi_\gamma(x^k)
			\]
			\hspace*{\algorithmicindent}
			for some \(H_k\in\hat\partial^2\varphi_\gamma(x^k)\)
		\State\label{step:FBN:LS}%
			Let \(\tau_k\) be the maximum in \(\set{2^{-i}}[i\in\N]\) such that
			\begin{equation}\label{eq:LS}
				\varphi_\gamma(x^{k+1})
				{}\leq{}
				\varphi_\gamma(x^k)
				{}-{}
				\sigma\|\Res(x^k)\|^2
			\end{equation}
			\hspace*{\algorithmicindent}
			where
			\(
				x^{k+1}
			{}\gets{}
				(1-\tau_k)\T(x^k)
				{}+{}
				\tau_k\bigl[x^k+d^k\bigr]
			\)
		\State
			\(k\gets k+1\) and go to \cref{step:FBN:1}
	\EndWhile
\end{algorithmic}

\end{algorithm}

In order to ease the computation of \(d^k\), we allow for inexact solutions of the linear system by introducing a tolerance \(\varepsilon_k>0\) and requiring
\(
	\|(H_k+\delta_k\I)d^k+\nabla\FBE(x^k)\|
{}\leq{}
	\varepsilon_k
\).
\begin{algorithm}[b]
	\algcaption{%
		({\sf CG}) Conjugate Gradient for computing the update direction%
	}%
	\label{alg:CG0}%
	\begin{algorithmic}[1]
	\Require
		\begin{tabular}[t]{@{}l@{}}
			\(\nabla\FBE(x^k)\);
			\(\delta_k\);
			\(\varepsilon_k\);
			\(d^{k-1}\) (set to \(0\) if \(k=0\))
			\\
			(generalized) directional derivatives
			\(\lambda\mapsto\dep{\prox_{\gamma g}}{\lambda}(\Fw{x^k})\)
			and
			\(\lambda\mapsto\dep{\nabla f}{\lambda}(x^k)\)
		\end{tabular}
	\Provide
		update direction \(d^k\)
	\Initialize
		\(e,p\gets-\nabla\FBE(x^k)\);~~
		warm start \(d^k\gets d^{k-1}\)
	\While{~\(\|e\|>\varepsilon_k\)~}
		\def\myVar#1{\fillwidthof[l]{d^k}{#1}}
		\State
			\(
				\myVar u
			{}\gets{}
				\dep{\nabla f}{p}(x^k)
			\)
		\State
			\(
				\myVar v
			{}\gets{}
				p-\gamma u
			\)
		\Comment{\(v=Q_\gamma(x^k)p\)}
		\State
			\(
				\myVar w
			{}\gets{}
				p-\dep{\prox_{\gamma g}}{v}(\Fw{x^k})
			\)
		\State
			\(
				\myVar z
			{}\gets{}
				\delta_kp
				{}+{}
				w-\gamma\dep{\nabla f}{w}(x^k)
			\)
		\Comment{\(z=H_kp\)}
		\State\label{step:pAp}\(
			\myVar\alpha
		{}\gets{}
			\|e\|^2 / \innprod pz
		\)
		\State\(
			\myVar{d^k}
		{}\gets{}
			d^k + \alpha p
		\),~~
		\(
			e^+
		{}\gets{}
			e - \alpha z
		\)
		\State\(
			\myVar p
		{}\gets{}
			e^+ + \bigl(\nicefrac{\|e^+\|}{\|e\|}\bigr)^2p
		\)
		\State\(
			\myVar e
		{}\gets{}
			e^+
		\)
	\EndWhile{}
\end{algorithmic}

\end{algorithm}
Since \(H_k+\delta_k\I\) is positive definite, inexact solutions of the linear system can be efficiently retrieved by means of \refCG, which only requires matrix-vector products and thus only (generalized) directional derivatives, namely, (generalized) derivatives (denoted as \(\dep{}{\lambda}\)) of the single-variable functions
\(
	t\mapsto\prox_{\gamma g}(x+t\lambda)
\)
and
\(
	t\mapsto\nabla f(x+t\lambda)
\),
as opposed to computing the full (generalized) Hessian matrix.
To further enhance computational efficiency, we may warm start the CG method with the previously computed direction, as eventually subsequent update directions are expected to have a small difference.
Notice that this warm starting does not ensure that the provided (inexact) solution \(d^k\) is a direction of descent for \(\FBE\); either way, this property is not required by the adopted linesearch, showing a considerable advantage over classical Armijo-type rules.
Putting all these facts together we obtain the proposed FBE-based truncated-Newton algorithm \refFBTN\ for convex composite minimization.
\begin{rem}[Adaptive variant when \(L_f\) is unknown]\label{rem:Adaptive}%
	In practice, no prior knowledge of the global Lipschitz constant \(L_f\) is required for \refFBTN[].
	In fact, replacing \(L_f\) with an initial estimate \(L>0\), the following instruction can be added at the beginning of each iteration, before \cref{step:FBN:1}:%
	\begin{algorithmic}[1]
	\makeatletter%
		\setcounter{algorithm}{1}%
		\setcounter{ALG@line}{-1}%
	\makeatother%
	\State
		\(\bar x^k\gets\T(x^k)\)
	\Statex
		{\bf while}~
		\(f(\bar x^k) > f(x^k)+\innprod{\nabla f(x^k)}{\bar x^k-x^k}+\frac L2\|\bar x^k-x^k\|^2\)
		~{\bf do}
		\item[]%
			\hspace*{\algorithmicindent}%
			\(\gamma\gets\nicefrac\gamma2\),~
			\(L\gets2L\),~
			\(\bar x^k\gets\T(x^k)\)
	\end{algorithmic}
	\stepcounter{algorithm}%
	Moreover, since positive definiteness of \(H_k+\delta_k\I\) is ensured only for \(\gamma\leq\nicefrac{1}{L_f}\) where \(L_f\) is the true Lipschitz constant of \(\nabla\FBE\) (cf. \cref{prop:PSDHess}), special care should be taken when applying \refCG[] in order to find the update direction \(d^k\).
	Specifically, \refCG[] should be stopped prematurely whenever \(\innprod pz\leq0\) in \cref{step:pAp}, in which case \(\gamma\gets\nicefrac\gamma2\), \(L\gets2L\) and the iteration should start again from \cref{step:FBN:1}.

	Whenever the quadratic bound \eqref{eq:LipBound} is violated with \(L\) in place of \(L_f\), the estimated Lipschitz constant \(L\) is increased, \(\gamma\) is decreased accordingly, and the proximal gradient point \(\bar x^k\) with the new stepsize \(\gamma\) is evaluated.
	Since replacing \(L_f\) with any \(L\geq L_f\) still satisfies \eqref{eq:LipBound}, it follows that \(L\) is incremented only a finite number of times.
	Therefore, there exists an iteration \(k_0\) starting from which \(\gamma\) and \(L\) are constant; in particular, all the convergence results here presented remain valid starting from iteration \(k_0\), at latest.
	Moreover, notice that this step does not increase the complexity of the algorithm, since both \(\bar x^k\) and \(\nabla f(x^k)\) are needed for the evaluation of \(\FBE(x^k)\).
\end{rem}

	\subsection{Subsequential and linear convergence}
		Before going through the convergence proofs let us spend a few lines to emphasize that \refFBTN[] is a well-defined scheme.
First, that a matrix \(H_k\) as in \cref{step:FBN:1} exists is due to the nonemptyness of \(\hat\partial^2\FBE(x^k)\) (cf. \cref{sec:2nd}).
Second, since \(\delta_k>0\) and \(H_k\succeq 0\) (cf. \cref{prop:PSDHess}) it follows that \(H_k+\delta_k\I\) is (symmetric and) positive definite, and thus \refCG[] is indeed applicable at \cref{step:FBN:CG}.

Having clarified this, the proof of the next result falls as a simplified version of \cite[Lem. 5.1 and Thm. 5.6]{themelis2018forward}; we elaborate on the details for the sake of self-inclusiveness.
To rule out trivialities, in the rest of the paper we consider the limiting case of infinite accuracy, that is \(\varepsilon=0\), and assume that the termination criterion \(\|\Res(x^k)\|=0\) is never met.
We shall also work under the assumption that a solution to the investigated problem \eqref{eq:GenProb} exists, thus in particular that the cost function \(\varphi\) is lower bounded.
\begin{thm}[Subsequential convergence]\label{thm:subseq}%
	Every accumulation point of the sequence \(\seq{x^k}\) generated by \refFBTN\ is optimal.
	\begin{proof}
		Observe that
		\[
			\FBE\bigl(x^k-\gamma\Res(x^k)\bigr)
		{}\overrel*[\leq]{\ref{thm:sandwich}}{}
			\FBE(x^k)
			{}-{}
			\gamma\tfrac{1-\gamma L_f}{2}\|\Res(x^k)\|^2
		{}<{}
			\FBE(x^k)
			{}-{}
			\sigma\|\Res(x^k)\|^2
		\]
		and that \(x^{k+1}\to\T(x^k)\) as \(\tau_k\to 0\).
		Continuity of \(\FBE\) ensures that for small enough \(\tau_k\) the linesearch condition \eqref{eq:LS} is satisfied, in fact, regardless of what \(d^k\) is.
		Therefore, for each \(k\) the stepsize \(\tau_k\) is decreased only a finite number of times.
		By telescoping the linesearch inequality \eqref{eq:LS} we obtain
		\[
			\sigma\sum_{k\in\N}{
				\|\Res(x^k)\|^2
			}
		{}\leq{}
			\sum_{k\in\N}\bigl[
				\varphi_\gamma(x^k)-\varphi_\gamma(x^{k+1})
			\bigr]
		{}\leq{}
			\varphi_\gamma(x^0)-\varphi_\star
		{}<{}
			\infty
		\]
		and in particular \(\Res(x^k)\to0\).
		Since \(\Res\) is continuous we infer that every accumulation point \(x_\star\) of \(\seq{x^k}\) satisfies \(\Res(x_\star)=0\), hence \(x_\star\in\argmin\varphi\), cf. \eqref{eq:R0}.
	\end{proof}
\end{thm}

\begin{rem}\label{rem:nonempty}%
	Since \refFBTN[] is a descent method on \(\FBE\), as ensured by the linesearch condition \eqref{eq:LS}, from \Cref{thm:LB} it follows that a sufficient condition for the existence of cluster points is having \(\varphi\) with bounded level sets or, equivalently, having \(\argmin\varphi\) bounded (cf. \cref{rem:coercive}).
\end{rem}

As a straightforward consequence of \Cref{thm:FBE_R2}, from the linesearch condition \eqref{eq:LS} we infer \(Q\)-linear decrease of the FBE along the iterates of \refFBTN[] provided that the original function \(\varphi\) has the quadratic growth property.
In particular, although the quadratic growth is a local property, \(Q\)-linear convergence holds globally, as described in the following result.

\begin{thm}[\(Q\)-linear convergence of {\refFBTN[]} under quadratic growth]\label{thm:Qlinear}%
	Suppose that \(\varphi\) satisfies the quadratic growth with constants \((\mu,\nu)\).
	Then, the iterates of \refFBTN\ decrease \(Q\)-linearly the value of \(\FBE\) as
	\[
		\FBE(x^{k+1})-\varphi_\star
	{}\leq{}
		\left(
			1
			{}-{}
			\tfrac{2\sigma\mu'}{
				\gamma\mu
				{}+{}
				2
				(2+\gamma\mu')
				(1+\gamma L_f)^2
			}
		\right)
		(\FBE(x^k)-\varphi_\star)
		\quad\forall k\in\N,
	\]
	where
	\[
		\mu'
	{}\coloneqq{}
		\begin{cases}[l @{~~} l]
			\mu
		&
			\text{if }\FBE(x_0)\leq\varphi_\star+\nu,
		\\
			\tfrac\mu2
			\min\set{
				1,\,
				\tfrac{\nu}{\FBE(x_0)-\varphi_\star-\nu}
			}
		&
			\text{otherwise.}
		\end{cases}
	\]
	\begin{proof}
		Since \refFBTN[] is a descent method on \(\FBE\), it holds that
		\(
			\seq{x^k}
		{}\subseteq{}
			\lev_{\leq\alpha}\FBE
		\)
		with \(\alpha=\FBE(x^0)\).
		It follows from \Cref{thm:QGglobal} that \(\varphi\) satisfies the quadratic growth condition with constants
		\(
			(\mu',\varphi_\gamma(x^0))
		\),
		with \(\mu'\) is as in the statement.
		The claim now follows from the inequality ensured by linesearch condition \eqref{eq:LS} combined with \Cref{thm:FBE_R2}.
	\end{proof}
\end{thm}

	\subsection{Superlinear convergence}
		In this section we provide sufficient conditions that enable superlinear convergence of \refFBTN[].
In the sequel, we will make use of the notion of \DEF{superlinear directions} that we define next.
\begin{defin}[Superlinear directions]
	Suppose that \(\mathcal X_\star\neq\emptyset\) and consider the iterates generated by \refFBTN.
	We say that \(\seq{d^k}\subset\R^n\) are \DEF{superlinearly convergent directions} if
	\begin{align*}
		\lim_{k\to\infty}{
			\frac{\dist(x^k+d^k,\mathcal X_\star)}{\dist(x^k,\mathcal X_\star)}
		}
	{}={} &
		0.
	\shortintertext{%
		If for some \(q>1\) the condition can be strengthened to
	}
		\limsup_{k\to\infty}{
			\frac{\dist(x^k+d^k,\mathcal X_\star)}{\dist(x^k,\mathcal X_\star)^q}
		}
	{}<{} &
		\infty
	\end{align*}
	then we say that \(\seq{d^k}\) are \DEF{superlinearly convergent directions with order \(q\)}.
\end{defin}
We remark that our definition of superlinear directions extends the one given in \cite[\S7.5]{facchinei2003finite} to cases in which \(\mathcal X_\star\) is not a singleton.
The next result consititutes a key component of the proposed methodology, as it shows that the proposed algorithm does not suffer from the \emph{Maratos' effect} \cite{maratos1978exact}, a well-known obstacle for fast local methods that inhibits the acceptance of the unit stepsize.
On the contrary, we will show that whenever the directions \(\seq{d^k}\) computed in \refFBTN[] are superlinear, then indeed the unit stepsize is eventually always accepted, and the algorithm reduces to a regularized version of the (undamped) linear Newton method \cite[Alg. 7.5.14]{facchinei2003finite} for solving \(\nabla\FBE(x)=0\) or, equivalently, \(\Res(x)=0\), and \(\dist(x^k,\mathcal X_\star)\) converges superlinearly.
\begin{thm}[Acceptance of the unit stepsize and superlinear convergence]\label{thm:tau1}%
	Consider the iterates generated by \refFBTN.
	Suppose that \(\varphi\) satisfies the quadratic growth (locally) and that \(\seq{d^k}\) are superlinearly convergent directions (with order \(q\)).
	Then, there exists \(\bar k\in\N\) such that
	\[
		\FBE(x^k+d^k)
	{}\leq{}
		\FBE(x^k)
		{}-{}
		\sigma\|\Res(x^k)\|^2
	\qquad
		\forall k\geq\bar k.
	\]
	In particular, eventually the iterates reduce to \(x^{k+1}=x^k+d^k\), and \(\dist(x^k,\mathcal X_\star)\) converges superlinearly (with order \(q\)).
	\begin{proof}
		Without loss of generality we may assume that \(\seq{x^k}\) and \(\seq{x^k+d^k}\) belong to a region in which quadratic growth holds.
		Denoting \(\varphi_\star\coloneqq\min\varphi\), since \(\FBE\) also satisfies the quadratic growth (cf. \cref{thm:QGequiv:phi}) if follows that
		\[
			\FBE(x^k)-\varphi_\star
		{}\geq{}
			\tfrac{\mu'}{2}\dist(x^k,\mathcal X_\star)^2
		\]
		for some constant \(\mu'>0\).
		Moreover, we know from \cref{thm:FBE_R2} that
		\[
			\FBE(x^k+d^k)-\varphi_\star
		{}\leq{}
			c\|\Res(x^k+d^k)\|^2
		{}\leq{}
			c'\dist(x^k+d^k,\mathcal X_\star)^2
		\]
		for some constants \(c,c'>0\), where in the second inequality we used Lipschitz continuity of \(\Res\) (\cref{thm:RLip}) together with the fact that \(\Res(x_\star)=0\) for all points \(x_\star\in \mathcal X_\star\).
		By combining the last two inequalities, we obtain
		\begin{equation}
			t_k
		{}\coloneqq{}
			\frac{
				\FBE(x^k+d^k)-\varphi_\star
			}{
				\FBE(x^k)-\varphi_\star
			}
		{}\leq{}
			\frac{
				2c'\dist(x^k+d^k,\mathcal X_\star)^2
			}{
				\mu'\dist(x^k,\mathcal X_\star)^2
			}
		{}\to{}
			0
		\quad
			\text{as }
			k\to\infty.
		\end{equation}
		Moreover,
		\begin{equation}
			\FBE(x^k)-\varphi_\star
		{}\geq{}
			\FBE(x^k)-\varphi(\T(x^k))
		{}\overrel*[\geq]{\ref{thm:geq}}{}
			\gamma\tfrac{1-\gamma L_f}{2}\|\Res(x^k)\|^2.
		\end{equation}
		Thus,
		\begin{qedalign*}
			\FBE(x^k+d^k)-\FBE(x^k)
		{}={} &
			\bigl[
				\FBE(x^k+d^k)-\varphi_\star
			\bigr]
			{}-{}
			\bigl[
				\FBE(x^k)-\varphi_\star
			\bigr]
		\\
		{}={} &
			(t_k-1)
			\bigl[
				\FBE(x^k)-\varphi_\star
			\bigr]
		\shortintertext{%
			and since \(t_k\to 0\), eventually it holds that
			\(
				t_k
			{}\leq{}
				1-\frac{2\sigma}{\gamma(1-\gamma L_f)}
			{}\in{}
				(0,1)
			\),
			resulting in
		}
		{}\leq{} &
			-\sigma\|\Res(x^k)\|^2.
		\end{qedalign*}
	\end{proof}
\end{thm}

\begin{thm}\label{thm:superlinear}%
	Consider the iterates generated by \refFBTN.
	Suppose that \(\varphi\) satisfies the quadratic growth (locally), and let \(x_\star\) be the limit point of \(\seq{x^k}\).\footnote{%
		As detailed in the proof, under the assumptions the limit point indeed exists.%
	}
	Then, \(\seq{d^k}\) are superlinearly convergent directions provided that
	\begin{enumerate}
	\item\label{thm:superlinear:QG}
		either \(\Res\) is strictly differentiable at \(x_\star\)\footnote{%
			From the chain rule of differentiation it follows that \(\Res\) is strictly differentiable at \(x_\star\) if \(\prox_{\gamma g}\) is strictly differentiable at \(\Fw{x_\star}\) (strict differentiability is closed under composition).%
		}
		and there exists \(D>0\) such that \(\|d^k\|\leq D\|\nabla\FBE(x^k)\|\) for all \(k\)'s,
	\item\label{thm:superlinear:strongmin}
		or \(\mathcal X_\star=\set{x_\star}\) and \(\prox_{\gamma g}\) is semismooth at \(\Fw{x_\star}\).
		In this case, if \(\prox_{\gamma g}\) is \(\vartheta\)-order semismooth at \(\Fw{x_\star}\) and \(\nabla^2f\) is \(\vartheta\)-H\"older continuous close to \(x_\star\), then the order of superlinear convergence is at least \(1+\min\set{\rho,\vartheta,\nu}\).
	\end{enumerate}
	\begin{proof}
		Due to \cref{thm:QGequiv,thm:Qlinear}, if \(\mathcal X_\star=\set{x_\star}\) then the sequence \(\seq{x^k}\) converges to \(x_\star\).
		Otherwise, the hypothesis ensure that
		\[
			\|x^{k+1}-x^k\|
		{}={}
			\tau_k\|d^k\|
		{}\leq{}
			D\|\nabla\FBE(x^k)\|
		{}\leq{}
			D\|\Res(x^k)\|,
		\]
		from which we infer that \(\seq{\|x^{k+1}-x^k\|}\) is \(R\)-linearly convergent, hence that \(\seq{x^k}\) is a Cauchy sequence, and again we conclude that the limit point \(x_\star\) indeed exists.
		Moreover, in light of \cref{thm:QGequiv} we have that \(\seq{x^k}\) is contained in a level set of \(\FBE\) where \(\FBE\) has quadratic growth.
		To establish a notation, let
		\(
			e^k
		{}\coloneqq{}
			[H_k+\delta_k\I]d^k
			{}+{}
			\nabla\FBE(x^k)
		\)
		be the error in solving the linear system at \cref{step:FBN:CG}, so that
		\begin{equation}
		\label{eq:e}
			\|e^k\|
		{}\leq{}
			\varepsilon_k
		{}\leq{}
			\|\nabla\FBE(x^k)\|^{1+\rho},
		\end{equation}
		(cf. \cref{step:FBN:1}), and let \(H_k=Q_\gamma(x^k)U_k\) for some \(U_k\in\mathcal\Res(x^k)\), see \eqref{eq:JR}.
		Let us now analyze the two cases separately.
		\begin{proofitemize}
		\item\ref{thm:superlinear:QG}~%
			Let \(x^k_\star\coloneqq\proj_{\mathcal X_\star}x^k\), so that \(\dist(x^k,\mathcal X_\star)=\|x^k-x^k_\star\|\).
			Recall that \(\nabla\FBE=Q_\gamma\Res\) and that \((1-\gamma L_f)\I\preceq Q_\gamma\preceq\I\).
			Since \(\Res(x^k_\star)=0\), from \cref{thm:RLip,thm:EB:Res} we infer that there exist \(r_1,r_2>0\) such that
			\begin{equation}
				\label{eq:Calmness}
					\|\Res(x^k)\|
				{}\geq{}
					r_1\dist(x^k,\mathcal X_\star)
			\quad\text{and}\quad
					\|\nabla\FBE(x^k)\|
				{}\leq{}
					r_2\dist(x^k,\mathcal X_\star).
			\end{equation}
			In particular, the assumption on \(d^k\) ensures that \(\|d^k\|=O\bigl(\dist(x^k,\mathcal X_\star)\bigr)\).
			We have
			\begin{align*}
				r_1\dist(x^k+d^k,\mathcal X_\star)
			{}\overrel*[\leq]{\eqref{eq:Calmness}}{} &
				\|\Res(x^k+d^k)\|
			\\
			{}\leq{} &
				\smashunderbrace{
					\|\Res(x^k+d^k)-\Res(x^k)-U_kd^k\|
				}{
					\text{(a)}
				}
				{}+{}
				\smashunderbrace{
					\|\Res(x^k)+U_kd^k\|
				}{
					\text{(b)}
				}.
			\end{align*}
			
			\vspace{.5\baselineskip}%
			As to quantity (a), we have
			\begin{align*}
				\text{(a)}
			{}\leq{} &
				\|\Res(x^k+d^k)-\Res(x^k)-J\Res(x_\star)d^k\|
				{}+{}
				\|U_k-JR(x_\star)\|
				\|d^k\|
			\\
			{}={} &
				o\bigl(\dist(x^k,\mathcal X_\star)\bigr),
			\shortintertext{%
				where we used strict differentiability and the fact that \(\partial_C\Res(x_\star)=\set{J\Res(x_\star)}\) \cite[Prop. 2.2.4]{clarke1990optimization} which implies \(U_k\to JR(x_\star)\).
				In order to bound (b), recall that \(\delta_k=\zeta\|\nabla\FBE(x^k)\|^\nu\) (cf. \cref{step:FBN:1}).
				Then,
			}
				\text{(b)}
			{}={} &
				\|Q_\gamma(x^k)^{-1}(e^k-\delta_kd^k)\|
			\\
			{}\overrel[\leq]{\eqref{eq:e}}{} &
				\tfrac{1}{1-\gamma L_f}
				\|\nabla\FBE(x^k)\|
				\left(
					\|\nabla\FBE(x^k)\|^\rho
					{}+{}
					\zeta\|\nabla\FBE(x^k)\|^{\nu-1}\|d^k\|
				\right).
			\\
			{}\overrel[\leq]{\eqref{eq:Calmness}}{} &
				\tfrac{r_2}{1-\gamma L_f}
				\dist(x^k,\mathcal X_\star)
				\left(
					r_2^\rho
					\dist(x^k,\mathcal X_\star)^\rho
					{}+{}
					\zeta\|\nabla\FBE(x^k)\|^{\nu-1}\|d^k\|
				\right)
			\\
			{}={} &
				O\bigl(\dist(x^k,\mathcal X_\star)^{1+\min\set{\rho,\nu}}\bigr),
			\end{align*}
			and we conclude that \(\dist(x^k+d^k,\mathcal X_\star)\leq\text{(a)}+\text{(b)}\leq o\bigl(\dist(x^k,\mathcal X_\star)\bigr)\).
		\item\ref{thm:superlinear:strongmin}~%
			In this case \(\dist(x^k,\mathcal X_\star)=\|x^k-x_\star\|\) and the assumption of (\(\vartheta\)-order) semismoothness ensures through \cref{thm:LNAR} that \(\mathcal\Res\) is a (\(\vartheta\)-order) LNA for \(\Res\) at \(x_\star\).
			Moreover, due to \cref{lem:sharpMin} there exists \(c>0\) such that \(\|[H_k+\delta_k\I]^{-1}\|\leq c\) for all \(k\)'s.
			We have
			\begin{align*}
				\|x^k+d^k-x_\star\|
			{}={} &
				\bigl\|
					x^k+[H_k+\delta_k\I]^{-1}(e^k-\nabla\FBE(x^k))-x_\star
				\bigr\|
			\\
			{}\leq{} &
				\bigl\|
					[H_k+\delta_k\I]^{-1}
				\bigr\|
				\bigl\|
					[H_k+\delta_k\I](x^k-x_\star)+e^k-\nabla\FBE(x^k)
				\bigr\|
			\\
				{}\leq{} &
					c
					\bigl\|
						H_k(x^k-x_\star)-\nabla\FBE(x^k)
					\bigr\|
					{}+{}
					c\delta_k
					\|x^k-x_\star\|
					{}+{}
					c\|e^k\|
			\\
				{}={} &
					c
					\bigl\|
						Q_\gamma(x^k)
						\bigl(
							\smashunderbracket{
								U_k(x^k-x_\star)-\Res(x^k)
							}{}
						\bigr)
					\bigr\|
					{}+{}
					c\delta_k
					\|x^k-x_\star\|
					{}+{}
					c\|e^k\|.
			\end{align*}
			Since \(\mathcal\Res\) is a LNA at \(x_\star\), it follows that the quantity emphasized in the bracket is a \(o(\|x^k-x_\star\|)\), whereas in case of a (\(\vartheta\)-order) LNA the tighter estimate \(O(\|x^k-x_\star\|^{1+\vartheta})\) holds.
			Combined with the fact that \(\delta_k=O(\|x^k-x_\star\|^\nu)\) and \(\|e^k\|=O(\|x^k-x_\star\|^{1+\rho})\), we conclude that \(\seq{d^k}\) are superlinearly convergent directions, and with order at least \(1+\min\set{\rho,\vartheta,\nu}\) in case of \(\vartheta\)-order semismoothness.
		\qedhere
		\end{proofitemize}
	\end{proof}
\end{thm}
Problems where the residual is (\(\vartheta\)-order) semismooth are quite common.
For instance, piecewise affine functions are everywhere strongly semismooth, as it is the case for the residual in lasso problems \cite{sopasakis2016accelerated}.
On the contrary, when the solution is not unique the condition \(\|d^k\|\leq D\|\nabla\FBE(x^k)\|\) (or, equivalently, \(\|d^k\|\leq D'\|\Res(x^k)\|\)) is trickier.
As detailed in \cite{zhou2005superlinear,zhou2006convergence}, this bound on the directions is ensured if \(\rho=1\) and for all iterates \(x^k\) and points \(x\) close enough to the limit point the following smoothness condition holds:
\begin{equation}\label{eq:Toh}
	\|\Res(x^k)+U_k(x-x^k)\|
{}\leq{}
	c\|x-x^k\|^2
\end{equation}
for some constant \(c>0\).
This condition is implied by and closely related to local Lipschitz differentiability of \(\Res\) and thus conservative.
We remark that, however, this can be weakened by requiring \(\rho\geq\nu\), and a notion of \(\vartheta\)-order semismoothness at the limit point with some degree of uniformity on the set of solutions \(\mathcal X_\star\), namely
\begin{equation}\label{eq:uLNA}
	\limsup_{
		\substack{
			x,x'\to x_\star
		\\
			x'\in\mathcal X_\star,\,
			x\neq x'
		\\
			U\in\mathcal\Res(x)
		}
	}{
		\frac{
			\|\Res(x)+U(x'-x)\|
		}{
			\|x'-x\|^{1+\vartheta}
		}
	}
{}<{}
	\infty
\end{equation}
for some \(\vartheta\in[\nu,1]\).
This weakened requirement comes from the observation that point \(x\) in \eqref{eq:Toh} is in fact \(x^k_\star\), the projection of \(x^k\) onto \(\mathcal X_\star\), set onto which \(\Res\) is constant (equal to \(0\)).
To see this, notice that \eqref{eq:uLNA} implies that
\(
	\|\Res(x^k)+U_k(x^k_\star-x^k)\|
{}\leq{}
	c\|x^k_\star-x^k\|^{1+\vartheta}
\)
for some \(c>0\).
In particular, mimicking the arguments in the cited references, since \(H_k\succeq0\) and \(\|Q_\gamma\|\leq1\), observe that
\begin{subequations}\label{subeq:norms}
	\begin{align}
	\label{eq:deltaInv}
		\|[H_k+\delta_k\I]^{-1}Q_\gamma(x^k)\|
	{}\leq{} &
		\|[H_k+\delta_k\I]^{-1}\|
	{}\leq{}
		\delta_k^{-1}
	\shortintertext{and}
	\label{eq:norm=2}
		\|[H_k+\delta_k\I]^{-1}H_k\|
	{}={} &
		\|\I-\delta_k[H_k+\delta_k\I]^{-1}\|
	{}\leq{}
		2.
	\end{align}
\end{subequations}
Therefore,
\begin{align*}
	\|d^k\|
{}={} &
	\bigl\|
		[H_k+\delta_k\I]^{-1}
		\bigl(
			e^k-\nabla\FBE(x^k)
		\bigr)
	\bigr\|
\\
{}\leq{} &
	\bigl\|
		[H_k+\delta_k\I]^{-1}
	\bigr\|
	\|e^k\|
	{}+{}
	\bigl\|
		[H_k+\delta_k\I]^{-1}
		Q_\gamma(x^k)
		\bigl(
			\Res(x^k)+U_k(x^k_\star-x^k)
		\bigr)
	\bigr\|
\\
&
	{}+{}
	\bigl\|
		[H_k+\delta_k\I]^{-1}
		H_k(x^k_\star-x^k)
	\bigr\|
\\
{}\overrel[\leq]{\eqref{subeq:norms}}{} &
	\delta_k^{-1}
	\|\nabla\FBE(x^k)\|^{1+\rho}
	{}+{}
	c\delta_k^{-1}
	\|x^k_\star-x^k\|^{1+\vartheta}
	{}+{}
	2\|x^k_\star-x^k\|
\\
{}={} &
	\zeta^{-1}
	\|\nabla\FBE(x^k)\|^{1+\rho-\nu}
	{}+{}
	c\zeta^{-1}
	\|x^k_\star-x^k\|^{1+\vartheta-\nu}
	{}+{}
	2\|x^k_\star-x^k\|
\\
{}={} &
	O\bigl(\dist(x^k,\mathcal X_\star)^{1+\min\set{0,\rho-\nu,\vartheta-\nu}}\bigr),
\end{align*}
which is indeed \(O(\|\nabla\FBE(x^k)\|)\) whenever \(\nu\leq\min\set{\vartheta,\rho}\).
Some comments are in order to expand on condition \eqref{eq:uLNA}.
\begin{enumerate}
\item
	If \(\mathcal X_\star=\set{\bar x}\) is a singleton, then \(x'\) is fixed to \(x_\star\) and the requirement reduces to \(\vartheta\)-order semismoothness at \(x_\star\).
\item
	This notion of uniformity is a \emph{local property}: for any \(\varepsilon>0\) the set \(\mathcal X_\star\) can be replaced by \(\mathcal X_\star\cap\ball{x_\star}{\varepsilon}\).
\item\label{thm:uLNA:B}%
	The condition \(U\in\mathcal\Res(x)\) in the limit can be replaced by
	\(
		U
	{}\in{}
		\hat{\mathcal\Res}(x)
	{}\coloneqq{}
		\set{\gamma^{-1}(\I-PQ_\gamma(x))}[
			P\in\partial_B\FB x
		]
	\),
	since \(\mathcal\Res(x)=\conv\bigl(\hat{\mathcal\Res}(x)\bigr)\).
\end{enumerate}
In particular, by exploiting this last condition it can be easily verified that if \(\Res\) is piecewise \(\vartheta\)-H\"older differentiable around \(x_\star\), then \eqref{eq:uLNA} holds, yet the stronger requirement \eqref{eq:Toh} in \cite{zhou2005superlinear,zhou2006convergence} does not.

	\section{Generalized Jacobians of proximal mappings}\label{sec:GenJac}
		In many interesting cases \(\prox_{\gamma g}\) is \(PC^1\) and thus semismooth.
\DEF{Piecewise quadratic} (PWQ) functions comprise a special but important class of convex functions whose proximal mapping is \(PC^1\).
A convex function \(g\) is called PWQ if \(\dom g\) can be represented as the union of finitely many polyhedral sets, relative to each of which \(g(x)\) is given by an expression of the form \(\frac 12\innprod{x}{Hx}+\innprod qx+c\) (\(H\in\R^{n\times n}\) must necessarily be symmetric positive semidefinite) \cite[Def. 10.20]{rockafellar2011variational}.
The class of PWQ functions is quite general since it includes \eg polyhedral norms, indicators and support functions of polyhedral sets, and it is closed under addition, composition with affine mappings, conjugation, inf-convolution and inf-projection \cite[Prop.s 10.22 and 11.32]{rockafellar2011variational}.
It turns out that the proximal mapping of a PWQ function is \emph{piecewise affine} (PWA) \cite[12.30]{rockafellar2011variational} (\(\R^n\) is partitioned in polyhedral sets relative to each of which \(\prox_{\gamma g}\) is an affine mapping), hence  strongly semismooth \cite[Prop. 7.4.7]{facchinei2003finite}.
Another example of a proximal mapping that is strongly semismooth is the projection operator over symmetric cones \cite{sun2002semismooth}.

A big class with semismooth proximal mapping is formed by the semi-algebraic functions.
We remind that a set \(A\subseteq\R^n\) is \DEF{semi-algebraic} if it can be expressed as
\[
	A
{}={}
	\bigcup_{i=1}^p{
		\bigcap_{j=1}^q\set{
			x\in\R^m
		}[
			P_{ij}(x)=0,\ Q_{ij}(x)<0
		]
	}
\]
for some polynomial functions \(\func{P_{ij},Q_{ij}}{\R^n}{\R}\), and that a function \(\func{h}{\R^n}{\Rinf{}^m}\) is \DEF{semi-algebraic} if \(\graph h\) is a semi-algebraic subset of \(\R^{n+m}\).
\begin{prop}
	If \(\func{g}{\R^n}{\Rinf}\) is semi-algebraic, then so are \(g^\gamma\) and \(\prox_{\gamma g}\).
	In particular, \(g^\gamma\) and \(\prox_{\gamma g}\) are semismooth.
	\begin{proof}
		Since \(g^\gamma\) and \(\prox_{\gamma g}\) are both Lipschitz continuous, semismoothness will follow once we show that they are semi-algebraic \cite[Rem. 4]{bolte2009tame}.
		Every polynomial is clearly semi-algebraic, and since the property is preserved under addition \cite[Prop. 2.2.6(ii)]{bochnak2013real}, the function \((x,w)\mapsto g(w)+\tfrac{1}{2\gamma}\|w-x\|^2\) is semi-algebraic.
		Moreover, since parametric minimization of a semi-algebraic function is still semi-algebraic (see \eg \cite[\S2]{attouch2013convergence}), it follows that the Moreau envelope \(g^\gamma\) is semi-algebraic and therefore so is
		\(
			h(x,w)
		{}\coloneqq{}
			g(w)+\tfrac{1}{2\gamma}\|w-x\|^2-g^\gamma(x)
		\).
		Notice that
		\(
			\prox_{\gamma g}(x)
		{}={}
			\set{w\in\R^n}[
				h(x,w)\leq 0
			]
		\),
		therefore
		\begin{align*}
			\graph\prox_{\gamma g}
		{}={} &
			\set{(x,\bar x)\in\R^n\times\R^n}[
				\prox_{\gamma g}(x)=\bar x
			]
		\\
		{}={} &
			\set{(x,\bar x)\in\R^n\times\R^n}[
				h(x,\bar x)\leq 0
			]
		\\
		{}={} &
			h^{-1}((-\infty,0])
		\end{align*}
		is a semi-algebraic set, since the interval \((-\infty,0]\) is clearly semi-algebraic and thus so is \(h^{-1}((-\infty,0])\) \cite[Prop. 2.2.7]{bochnak2013real}.
	\end{proof}
\end{prop}
In fact, with the same arguments it can be shown that the result still holds if `semi-algebraic' is replaced with the broader notion of `\DEF{tame}', see \cite{bolte2009tame}.
Other conditions that guarantee semismoothness of the proximal mapping can be found in \cite{meng2005semismoothness,meng2008lagrangian,meng2008moreau,mifflin1999properties}.
The rest of the section is devoted to collecting explicit formulas of \(\partial_C\prox_{\gamma g}\) for many known useful instances of convex functions \(g\).

		\subsection{Properties}
			\begin{proxexamples}[itemsep=5pt]
\item[Separable functions]\label{es:Separ}%
	
	Whenever \(g\) is (block) separable, \ie \(g(x)=\sum_{i=1}^N g_i(x_i)\), \(x_i\in\R^{n_i}\), \(\sum_{i=1}^N n_i=n\), then every \(P\in\partial_C(\prox_{\gamma g})(x)\) is a (block) diagonal matrix.
	This has favorable computational implications especially for large-scale problems.
	For example, if \(g\) is the \(\ell_1\) norm or the indicator function of a box, then the elements of \(\partial_C\prox_{\gamma g}(x)\) (or \(\partial_B\prox_{\gamma g}(x)\)) are diagonal matrices with diagonal elements in \([0,1]\) (or in \(\set{0,1}\)).
\item[Convex conjugate]\label{es:MorDec}

	With a simple application of the Moreau's decomposition \cite[Thm. 14.3(ii)]{bauschke2017convex}, all elements of \(\partial_C\prox_{\gamma\conj g}\) are readily available as long as one can compute \(\partial_C\prox_{\nicefrac g\gamma}\).
	Specifically,
	\[
		\partial_C(\prox_{\gamma\conj g})(x)
	{}={}
		\I-\partial_C(\prox_{\nicefrac g\gamma})(\nicefrac x\gamma).
	\]
\item[Support function]

	The \DEF{support function} of a nonempty closed and convex set \(D\) is the proper convex and lsc function \(\sigma_D(x)\coloneqq\sup_{y\in D}\innprod xy\).
	Alternatively, \(\sigma_D\) can be expressed as the convex conjugate of the indicator function \(\delta_D\), and one can use the results of \Cref{es:MorDec} to find that
	\[
		\partial_C(\prox_{\gamma g})(x)
	{}={}
		\I-\partial_C(\proj_D)(\nicefrac x\gamma).
	\]
	\Cref{ex:IndFun} offers a rich list of sets \(D\) for which a close form expression exists.
\item[Spectral functions]
	
	The eigenvalue function \(\func{\lambda}{\sym(\R^{n\times n})}{\R^n}\) returns the vector of eigenvalues of a symmetric matrix in nonincreasing order.
	\emph{Spectral functions} are of the form
	\begin{equation}\label{eq:SpecFun}
		\func{G\coloneqq h\circ\lambda}{\sym(\R^{n\times n})}{\Rinf}.
	\end{equation}
	where \(\func h{\R^n}{\Rinf}\) is proper, lsc, convex and \emph{symmetric}, \ie invariant under coordinate permutations \cite{lewis1996convex}.
	Such \(G\) inherits most of the properties of \(h\) \cite{lewis1996derivatives,lewis2001twice}; in particular, its proximal mapping is \cite[\S6.7]{parikh2014proximal}
	\[
		\prox_{\gamma G}(X)
	{}={}
		Q\diag(\prox_{\gamma h}(\lambda(X)))\trans Q,
	\]
	where \(X=Q\diag(\lambda(X))\trans Q\) is the spectral decomposition of \(X\) (\(Q\) is an orthogonal matrix).
	If, additionally,
	\begin{equation}\label{eq:SymSep}
		h(x)=g(x_1)+\cdots+g(x_N)
	\end{equation}
	for some \(\func g\R\Rinf\),
	then
	\[
		\prox_{\gamma h}(x)
	{}={}
		(\prox_{\gamma g}(x_1),\ldots,\prox_{\gamma g}(x_N)),
	\] 
	and therefore the proximal mapping of \(G\) can be expressed as
	\begin{equation}\label{eq:ProxSpec}
		\prox_{\gamma G}(X)
	{}={}
		Q\diag(\prox_{\gamma g}(\lambda_1(X)),\ldots,\prox_{\gamma g}(\lambda_n(X)))\trans Q,
	\end{equation}
	\cite[Chap. V]{bhatia1997matrix}, \cite[Sec. 6.2]{horn1994topics}.
	Now we can use the theory of nonsmooth symmetric matrix-valued functions developed in \cite{chen2003analysis} to analyze differentiability properties of \(\prox_{\gamma G}\).
	In particular, \(\prox_{\gamma G}\) is (strongly) semismooth at \(X\) iff \(\prox_{\gamma g}\) is (strongly) semismooth at the eigenvalues of \(X\) \cite[Prop. 4.10]{chen2003analysis}.
	Moreover, for any \(X\in\sym(\R^{n\times n})\) and \(P\in\partial_B(\prox_{\gamma G})(X)\) we have \cite[Lem. 4.7]{chen2003analysis} 
	\begin{equation}\label{eq:JacSpec}
		P(X)
	{}={}
		Q
		\left(
			\Omega^{\gamma g}_{\lambda,\lambda}\odot(\trans QXQ)
		\right)
		\trans Q,
	\end{equation}
	where \(\odot\) denotes the Hadamard product and for vectors \(u,v\in\R^n\) we defined \(\Omega^{\gamma g}_{u,v}\) as the \(n\times n\) matrix
	\begin{equation}\label{eq:GammaJac}
		(\Omega^{\gamma g}_{u,v})_{ij}
	{}\coloneqq{}
		\begin{cases}[l @{~~} l]
			\partial_B\prox_{\gamma g}(u_i)
			&
			\text{if } u_i=v_j,
		\\[4pt]
			\set{
				\frac{
					\prox_{\gamma g}(u_i)-\prox_{\gamma g}(v_j)
				}{
					u_i-v_j
				}
			}
			&
			\text{otherwise.}
		\end{cases}
	\end{equation}
\item[Orthogonally invariant functions]
	A function \(\func{G}{\R^{m\times n}}{\Rinf}\) is called \emph{orthogonally invariant} if
	\(
		G(UX\trans V)=G(X)
	\)
	for all \(X\in\R^{m\times n}\) and orthogonal matrices \(U\in\R^{m\times m}\), \(V\in\R^{n\times n}\).\footnote{%
		In case of complex-valued matrices, functions of this form are known as \emph{unitarily invariant} \cite{lewis1995convex}.
	}

	A function \(\func{h}{\R^q}{\Rinf}\) is \emph{absolutely symmetric} if \(h(Qx)=h(x)\) for all \(x\in\R^q\) and any generalized permutation matrix \(Q\), \ie a matrix \(Q\in\R^{q\times q}\) that has exactly one nonzero entry in each row and each column, that entry being \(\pm 1\) \cite{lewis1995convex}.
	There is a one-to-one correspondence between orthogonally invariant functions on \(\R^{m\times n}\) and absolutely symmetric functions on \(\R^q\).
	Specifically, if \(G\) is orthogonally invariant then 
	\[
		G(X)=h(\sigma(X))
	\]
	for the absolutely symmetric function \(h(x)=G(\diag(x))\).
	Here, for \(X\in\R^{m\times n}\) and \(q\coloneqq\min\set{m,n}\) the spectral function \(\func{\sigma}{\R^{m\times n}}{\R^q}\) returns the vector of its singular values in nonincreasing order.
	Conversely, if \(h\) is absolutely symmetric then \(G(X)=h(\sigma(X))\) is orthogonally invariant.
	Therefore, convex analytic and generalized differentiability properties of orthogonally invariant functions can be easily derived from those of the corresponding absolutely symmetric functions \cite{lewis1995convex}.
	For example, assuming for simplicity that \(m\leq n\), the proximal mapping of \(G\) is given by \cite[Sec. 6.7]{parikh2014proximal}
	\[
		\prox_{\gamma G}(X)
	{}={}
		U\diag(\prox_{\gamma h}(\sigma(X)))\trans{V_1},
	\]
	where
	\(
		X
	{}={}
		U
		\begin{bmatrix}\diag(\sigma(X))&0\end{bmatrix}
		\trans{\begin{bmatrix}V_1 & V_2\end{bmatrix}}
	\)
	is the singular value decomposition of \(X\).
	If we further assume that \(h\) has a separable form as in \eqref{eq:SymSep}, then
	\begin{equation}\label{eq:ProxSpec2}
		\prox_{\gamma G}(X)=U\Sigma_g(X)\trans{V_1},
	\end{equation}
	where \(\Sigma_g(X)=\diag(\prox_{\gamma g}(\sigma_1(X)),\ldots,\prox_{\gamma g}(\sigma_n(X)))\).
	Functions of this form are called \emph{nonsymmetric matrix-valued functions}.
	We also assume that \(g\) is a non-negative function such that \(g(0)=0\).
	This implies that \(\prox_{\gamma g}(0)=0\) and guarantees that the nonsymmetric matrix-valued function \eqref{eq:ProxSpec2} is well defined \cite[Prop. 2.1.1]{yang2009study}.
	Now we can use the results of \cite[\S2]{yang2009study} to draw conclusions about generalized differentiability properties of \(\prox_{\gamma G}\).
	For example, through \cite[Thm. 2.27]{yang2009study} we have that \(\prox_{\gamma G}\) is continuously differentiable at \(X\) if and only if \(\prox_{\gamma g}\) is continuously differentiable at the singular values of \(X\).
	Furthermore, \(\prox_{\gamma G}\) is (strongly) semismooth at \(X\) if \(\prox_{\gamma g}\) is (strongly) semismooth at the singular values of \(X\) \cite[Thm. 2.3.11]{yang2009study}. 

	For any \(X\in\R^{m\times n}\) the generalized Jacobian \(\partial_B(\prox_{\gamma G})(X)\) is well defined and nonempty, and any \(P\in\partial_B(\prox_{\gamma G})(X)\) acts on \(H\in\R^{m\times n}\) as \cite[Prop. 2.3.7]{yang2009study}
%
%
	\begin{equation}\label{eq:JacNS}
		P(H)
	{}={}
		U
		\begin{bmatrix}
			\left(
				\Omega^{\gamma g}_{\sigma,\sigma}\odot\left(\frac{H_1+\trans{H_1}}{2}\right)
				{}+{}
				\Omega^{\gamma g}_{\sigma,-\sigma}\odot\left(\frac{H_1-\trans{H_1}}{2}\right)
			\right),
			&
			(\Omega^{\gamma g}_{\sigma,0}\odot H_2)
		\end{bmatrix}
		\trans V,
	\end{equation}
	where \(V=[V_1~V_2]\), \(H_1=\trans UHV_1\in\R^{m\times m}\), \(H_2=\trans UHV_2\in\R^{m\times(n-m)}\) and matrices \(\Omega\) are as in \eqref{eq:GammaJac}.
\end{proxexamples}

		\subsection{Indicator functions}\label{ex:IndFun}%
			Smooth constrained convex problems
\begin{equation}\label{eq:CP}
	\minimize_{x\in\R^n} f(x)
\quad\stt{}
	x\in D
\end{equation}
can be cast in the composite form \eqref{eq:GenProb} by encoding the feasible set \(D\) with the indicator function \(g=\indicator_D\).
Whenever \(\proj_D\) is efficiently computable, then algorithms like the forward-backward splitting \eqref{eq:FBS} can be conveniently considered.
In the following we analyze the generalized Jacobian of some of such projections.

\begin{proxexamples}
\item[Affine sets]\label{ex:aff}
	\(D=\set{x\in\R^n}[Ax=b]\) for some \(A\in\R^{m\times n}\) and \(b\in\R^m\).

	In this case, \(\proj_D(x)=x-\pinv A(Ax-b)\) where \(\pinv A\) is the Moore-Penrose pseudoinverse of \(A\).
	For example, if \(A\) is surjective (\ie, it has full row rank and thus \(m\leq n\)), then \(\pinv A=\trans A(A\trans A)^{-1}\), whereas if it is injective (\ie, it has full column rank and thus \(m\geq n\)), then \(\pinv A=(\trans AA)^{-1}\trans A\).
	Obviously \(\proj_D\) is an affine mapping, thus everywhere differentiable with
	\[
		\partial_C(\proj_D)(x)
	{}={}
		\partial_B(\proj_D)(x)
	{}={}
		\set{\nabla\proj_D(x)}
	{}={}
		\set{\I-\pinv AA}.
	\]
\item[Polyhedral sets]\label{ex:poly}
	\(D=\set{x\in\R^n}[Ax=b,\ Cx\leq d]\), for some \(A\in\R^{p\times n}\), \(b\in\R^p\), \(C\in\R^{q\times n}\) and \(d\in\R^q\).

	It is well known that \(\proj_D\) is piecewise affine.
	In particular, let 
	\[
		\mathcal I_D
	{}\coloneqq{}
		\set{I\subseteq\set{1\ldots q}}[
			\exists x\in\R^n
			{}:{}
			Ax=b,\ 
			C_{i\cdot}x=d_i\ i\in I,\ 
			C_{j\cdot}x<d_j\ j\notin I
		].
	\]
	Then, the \DEF{faces} of \(D\) can be indexed with the elements of \(\mathcal I\) \cite[Prop. 2.1.3]{scholtes2012piecewise}: for each \(I\in\mathcal I_D\) let 
	\begin{align*}
		F_I
	{}\coloneqq{} &
		\set{x\in D}[C_{i\cdot}x=d_i,\ i\in I]
	\shortintertext{%
		be the \(I\)-th face of \(D\),
	}
		S_I
	{}\coloneqq{} &
		\aff F_I
	{}={}
		\set{x\in\R^n}[Ax=b,\ C_{i\cdot}x=d_i,\ i\in I]
	\shortintertext{%
		be the hyperplane containing the \(I\)-th face of \(D\),
	}
		N_I
	{}\coloneqq{} &\textstyle
		\range\trans A
		{}+{}
		\cone\set{
			\trans{C_{I\cdot}}
		}
	\shortintertext{%
		be the normal cone to any point in the relative interior of \(F_I\) \cite[Eq. (2.44)]{scholtes2012piecewise},\footnotemark{} and
	}
		R_I
	{}\coloneqq{} &
		F_I+N_I.
	\end{align*}
	\footnotetext{%
		Consistently with the definition in \cite{scholtes2012piecewise}, the polyhedron \(P\) can equivalently be expressed by means of only inequalities as
		\(
			P
		{}={}
			\set{x\in\R^n}[Ax\leq b,~-Ax\leq -b,~Cx\leq b]
		\),
		resulting indeed in
		\(
			\cone [\trans A,~-\trans A,~\trans C]
		{}={}
			\range\trans A+\cone\trans C
		\).
	}%
	We then have \(\proj_D(x)\in\set{\proj_{S_I}(x)}[I\in\mathcal I_D]\), \ie \(\proj_D\) is a piecewise affine function.
	The affine  pieces of \(\proj_D\) are the projections on the corresponding affine subspaces \(S_I\) (cf. \Cref{ex:aff}).
	In fact, for each \(x\in R_I\) we have \(\proj_D(x)=\proj_{S_I}(x)\), each \(R_I\) is full dimensional and \(\R^n=\bigcup_{I\in\mathcal I_D}R_I\) \cite[Prop.s 2.4.4 and 2.4.5]{scholtes2012piecewise}.
	For each \(I\in \mathcal I_D\) let
	\begin{equation}\label{eq:PI}
		P_I
	{}\coloneqq{}
		\nabla\proj_{S_I}
	{}={}
		\I-\pinv{\binom{A}{C_I}}\binom{A}{C_I},
	\end{equation}
	and for each \(x\in\R^n\) let
	\[
		\mathcal I_D(x)
	{}\coloneqq{}
		\set{I\in\mathcal I_D}[x\in R_I].
	\]
	Then,
	\[
		\partial_C(\proj_D)(x)
	{}={}
		\conv\partial_B(\proj_D)(x)
	{}={}
		\conv\set{P_I}[I\in \mathcal I_D(x)].
	\]
	Therefore, an element of \(\partial_B\proj_D(x)\) is \(P_I\) as in \eqref{eq:PI} where
	\(
		I
	{}={}
		\set{i}[
			C_i\bar x=d_i
		]
	\)
	is the set of active constraints of \(\bar x=\proj_D(x)\).
	For a more general analysis we refer the reader to \cite{han1997newton,li2017efficient}.
\item[Halfspaces]
	\(H=\set{x\in\R^n}[\innprod ax\leq b]\) for some \(a\in\R^n\) and \(b\in\R\).
	
	Then, denoting the \DEF{positive part} of \(r\in\R\) as \([r]_+\coloneqq\max\set{0,r}\),
	\begin{align*}
		\proj_H(x)
	{}={} &
		x-\tfrac{[\innprod ax-b]_+}{\|a\|^2}a
	\shortintertext{and}
		\partial_C(\proj_H)(x)
	{}={} &
		\begin{cases}[l @{~~} l]
			\set{\I-\|a\|^{-2}a\trans a} & \text{if } x\notin H,
		\\
			\set{\I} & \text{if } \innprod ax<b,
		\\
			\conv\set{\I,\I-\|a\|^{-2}a\trans a} & \text{if } \innprod ax=b.
		\end{cases}
	\end{align*}
\item[Boxes]\label{ex:projBox}
	\(D=\set{x\in\R^n}[\ell\leq x\leq u]\) for some \(\ell,u\in[-\infty,\infty]^n\).
	
	We have
	\begin{align*}
		\proj_D(x)
	{}={} &
		\min\set{\max\set{x,\ell},u},
	\shortintertext{%
		and since the corresponding indicator function \(\delta_D\) is separable, every element of \(\partial_C(\proj_D)(x)\) is diagonal with (cf. \cref{es:Separ})
	}
		\partial_C(\proj_D)(x)_{ii}
	{}={} &
		\begin{cases}[l @{~~} l]
			\set0 & \text{if } x\notin D,\\
			\set1 & \text{if } \ell_i<x_i<u_i,\\{}
			[0,1] & \text{if } x_i\in\set{\ell_i,u_i}.
		\end{cases}
	\end{align*}
\item[Unit simplex]\label{ex:projSimplex}%
	\(D=\set{x\in\R^n}[x\geq 0,\ \sum_{i=1}^nx_i=1]\).
	
	By writing down the optimality conditions for the corresponding projection problem, one can easily see that
	\[
		\proj_D(x)=[x-\lambda{\bf 1}]_+,
	\]
	where \(\lambda\) solves \(\innprod{\mathbf 1}{[x-\lambda\mathbf 1]_+}=1\).
	Since the unit simplex is a polyhedral set, we are dealing with a special case of \Cref{ex:poly}, where \(A=\trans{\bf 1}\), \(b=1\), \(C=-\I\) and \(d=0\).
	Therefore, in order to calculate an element of the generalized Jacobian of the projection, we first compute \(\proj_D(x)\) and then determine the set of active indices \(J\coloneqq\set{i}[\proj_D(x)_i=0]\).
	An element \(P\in\partial_B(\proj_D)(x)\) is given by
	\[
		P_{ij}
	{}={}
		\begin{cases}[l @{~~} l]
			\delta_{i,j}-\frac{1}{n-|J|} & \text{if } i,j\notin J,
		\\
			0 & \text{otherwise,}
		\end{cases}
	\]
	where \(|J|\) denotes the cardinality of the set \(J\).
	Notice that \(P\) is block diagonal after a permutation of rows and columns. 
\item[Euclidean unit ball]
	\(B=\cball{0}{1}\).
	
	We have
	\begin{align*}
		\proj_B(x)
	{}={} &
		\begin{cases}[l @{~~} l]
			x & \text{if } x\in B,\\
			\nicefrac{x}{\|x\|} & \text{otherwise,}
		\end{cases}
	\shortintertext{and}
		\partial_C(\proj_B)(x)
	{}={} &
		\begin{cases}[l @{~~} l]
			\set{\I} & \text{if } \|x\|<1,
		\\
		\conv\set{\|x\|^{-1}(\I-w\trans w),\,\I} & \text{if } \|x\|=1,
		\\
			\set{\|x\|^{-1}(\I-w\trans w)} & \text{if } x\notin B,
		\end{cases}
	\end{align*}
	where \(w\coloneqq\nicefrac{x}{\|x\|}\).
\item[Second-order cone]
	\(
		\mathcal K
	{}={}
		\set{(x_0,\bar x)\in\R\times\R^{n-1}}[x_0\geq\|\bar x\|]
	\).
	
	Let \(x\coloneqq(x_0,\bar x)\), and for \(w\in\R^n\) and \(\alpha\in\R\) define
	\[
		M_{w,\alpha}
	{}\coloneqq{}
		\tfrac12\begin{bmatrix}
			1 & \trans w\\
			w & (1-\alpha)\I_{n-1}+\alpha w\trans w
		\end{bmatrix}.
	\]
	Then,
	\(
		\partial_C(\proj_{\mathcal K})(x)
	{}={}
		\conv(\partial_B(\proj_{\mathcal K})(x))
	\)
	where, for \(\bar w\coloneqq\nicefrac{\bar x}{\|\bar x\|}\) and \(\bar\alpha\coloneqq-\nicefrac{x_0}{\|\bar x\|}\), we have \cite[Lem. 2.6]{kanzow2009local}
	\[
		\partial_B(\proj_{\mathcal K})(x)
	{}={}
		\begin{cases}[l@{~~\text{if }}l]
			\set0
			&
			x_0<-\|\bar x\|,
		\\
			\set{\I_n}
			&
			x_0>\|\bar x\|,
		\\
			\set{M_{\bar w,\bar\alpha}}
			&
			-\|\bar x\|<x_0<\|\bar x\|,
		\\
			\set{\I_n,\,M_{\bar w,\bar\alpha}}
			&
			x_0=\|\bar x\|\neq 0,
		\\
			\set{0,\,M_{\bar w,\bar\alpha}}
			&
			x_0=-\|\bar x\|\neq 0,
		\\
			\set{0,\I_n}
			{}\cup{}
			\set{M_{w,\alpha}}[
				|\alpha|\leq1,\,\|w\|\leq 1
			]
			&
			x_0=\bar x= 0.
		\end{cases}
	\]
\item[Positive semidefinite cone]
	\(\mathcal S_+=\sym_+(\R^{n\times n})\).
	
	For any symmetric matrix \(M\) it holds that
	\[
		\proj_{\mathcal S_+}(M)
	{}={}
		Q[\diag(\lambda)]_+\trans Q,
	\]
	where \(M=Q\diag(\lambda)\trans Q\) is any spectral decomposition of \(M\).
	This coincides with \eqref{eq:ProxSpec}, as \(\indicator_{\mathcal S_+}\) can be expressed as in \eqref{eq:SpecFun}, where \(h\) has the separable form \eqref{eq:SymSep} with \(g=\delta_{\R_+}\), so that for \(r\in\R\) we have
	\[
		\prox_{\gamma g}(r)=[r]_+
	\quad\text{and}\quad
		\partial_B(\prox_{\gamma g})(r)
	{}={}
		\begin{cases}[l@{~~}l]
			\set 0 & \text{if } r<0,
		\\
			\set{0,1} & \text{if } r=0,
		\\
			\set 1 & \text{if } r>0.
		\end{cases}
	\]
	An element of \(\partial_B\proj_{\sym_+(\R^{n\times n})}(X)\) is thus given by \eqref{eq:JacSpec}.
\end{proxexamples}

		\subsection{Norms}
			\begin{proxexamples}
\item[\texorpdfstring{\(\ell_1\)}{l1} norm]\label{ex:EllOne}
	\(g(x)=\|x\|_1\).
	
	The proximal mapping is the well known soft-thresholding operator
	\[
		(\prox_{\gamma g}(x))_i
	{}={}
		\sign(x_i)[|x_i|-\gamma]_+,
	\quad
		i=1,\ldots,n.
	\]
	Function \(g\) is separable, and thus every element of \(\partial_B(\prox_{\gamma g})\) is a diagonal matrix, cf. \cref{es:Separ}.
	Specifically, the nonzero elements are
	\[
		\partial_C(\prox_{\gamma g})(x)_{ii}
	{}={}
		\begin{cases}[l @{~~} l]
			\set1 & \text{if } |x_i|>\gamma,\\{}
			[0,1] & \text{if } |x_i|=\gamma,\\
			\set0 & \text{if } |x_i|<\gamma.
		\end{cases}
	\]
	We could also arrive to the same conclusion by applying the Moreau decomposition of \Cref{es:MorDec} to the function of \Cref{ex:projBox} with \(u=-\ell={\bf 1}_n\), since the \(\ell_1\) norm is the conjugate of the indicator of the \(\ell_\infty\)-norm ball.
\item[\texorpdfstring{\(\ell_\infty\)}{l-infinity} norm]
	\(g(x)=\|x\|_\infty\).

	Function \(g\) is the convex conjugate of the indicator of the unit simplex \(D\) analyzed in \Cref{ex:projSimplex}.
	From the Moreau decomposition, see \cref{es:MorDec}, we obtain
	\[
		\partial_C(\prox_{\gamma g})(x)
	{}={}
		\I-\partial_C(\proj_D)(\nicefrac x\gamma).
	\]
	Then,
	\(
		\proj_D(\nicefrac x\gamma)
	{}={}
		[\nicefrac x\gamma-\lambda\mathbf{1}]_+
	\)
	where \(\lambda\in\R\) solves \(\innprod{\mathbf{1}}{[\nicefrac x\gamma-\lambda\mathbf{1}]_+}=1\).
	Let \(J=\set{i}[\proj_D(\nicefrac x\gamma)_i=0]\), then an element of \(\partial_B(\prox_{\gamma g})(x)\) is given by
	\[
		P_{ij}
	{}={}		
		\begin{cases}[l @{~~} l]
			\frac{1}{n-|J|} & \text{if } i,j\notin J,
		\\
			\delta_{i,j} & \text{otherwise.}
		\end{cases}
	\]
\item[Euclidean norm]
	\(g(x)=\|x\|\).
	
	The proximal mapping is given by
	\begin{equation}\label{eq:SoftThresholding}
		\prox_{\gamma g}(x)
	{}={}
		\begin{cases}[l @{~~} l]
			(1-\gamma\|x\|^{-1})x & \text{if } \|x\|\geq\gamma,\\
			0 & \text{otherwise.}
		\end{cases}
	\end{equation}
	Since \(\prox_{\gamma g}\) is a \(PC^1\) mapping, its \(B\)-subdifferential can be computed by simply computing the Jacobians of its smooth pieces. Specifically, denoting \(w=\nicefrac{x}{\|x\|}\) we have
	\[
		\partial_C(\prox_{\gamma g})(x)
	{}={}
		\begin{cases}[l @{~~} l]
			\set{\I-\gamma\|x\|^{-1}(\I-w\trans w)} & \text{if } \|x\|>\gamma,
		\\
			\set{0} & \text{if } \|x\|<\gamma,
		\\
			\conv\set{\I-\gamma\|x\|^{-1}(\I-w\trans w),0} & \text{otherwise.}
		\end{cases}
	\]
\item[Sum of Euclidean norms]
	\(g(x)=\sum_{s\in\mathcal S}\|x_s\|\), where \(\mathcal S\) is a partition of \(\set{1,\ldots,n}\).
	
	Differently from the \(\ell_1\)-norm which induces sparsity on the whole vector, this function serves as regularizer to induce group sparsity \cite{yuan2006model}.
	For \(s\in\mathcal S\), the components of the proximal mapping indexed by \(s\) are
	\[
		(\prox_{\gamma g}(x))_s
	{}={}
		(1-\gamma\|x_s\|^{-1})_+x_s.
	\]
	Any \(P\in\partial_B(\prox_{\gamma g})(x)\) is block diagonal with the \(s\)-block equal to
	\[
		P_s
	{}={}
		\begin{cases}[l @{~~} l]
			\I-\gamma\|x_s\|^{-1}\bigl(\I-\|x_s\|^{-2}x_s\trans{x_s})
			&
			\text{if }\|x_s\|>\gamma,
		\\
			\I
			&
			\text{if }\|x_s\|<\gamma,
		\\
			\text{any of these two matrices}
			&
			\text{if }\|x_s\|=\gamma.
		\end{cases}
	\]
\item[Matrix nuclear norm] \(G(X)=\|X\|_\star\) for \(X\in\R^{m\times n}\).

	The \emph{nuclear norm} returns the sum of the singular values of a matrix \(X\in\R^{m\times n}\), \ie \(G(X)=\sum_{i=1}^m\sigma_i(X)\) (for simplicity we are assuming that \(m\leq n\)).
	It serves as a convex surrogate for the rank, and has found many applications in systems and control theory, including system identification and model reduction \cite{fazel2001rank,fazel2002matrix,fazel2004rank,liu2010interior,recht2010guaranteed}.
	Other fields of application include \emph{matrix completion problems} arising in machine learning \cite{srebro2004learning,rennie2005fast} and computer vision \cite{tomasi1992shape,morita1997sequential}, and \emph{nonnegative matrix factorization problems} arising in data mining \cite{elden2007matrix}.

	The nuclear norm can be expressed as \(G(X)=h(\sigma(X))\), where \(h(x)=\|x\|_1\) is absolutely symmetric and separable.
	Specifically, it takes the form \eqref{eq:SymSep} with \(g=|{}\cdot{}|\), for which \(g(0)=0\) and \(0\in\partial g(0)\), and whose proximal mapping is the soft-thresholding operator.
	In fact, since the case of interest here is \(x\geq 0\) (because \(\sigma_i(X)\geq 0\)), we have \(\prox_{\gamma g}(x)=[x-\gamma]_+\), cf. \eqref{eq:SoftThresholding}.
	Consequently, the proximal mapping of \(\|X\|_\star\) is given by \eqref{eq:ProxSpec2} with
	\[
		\Sigma_g(X)=\diag([\sigma_1(X)-\gamma]_+,\ldots,[\sigma_m(X)-\gamma]_+).
	\]
	For \(x\in\R_+\) we have that 
	\begin{equation}\label{eq:subSoft}
		\partial_C(\prox_{\gamma g})(x)
	{}={}
		\begin{cases}[l @{~~} l]
			0     & \text{if } 0\leq x<\gamma,\\{}
			[0,1] & \text{if } x=\gamma,\\
			1     & \text{if } x>\gamma,
		\end{cases}
	\end{equation}
	then \(\partial_B(\prox_{\gamma G})(X)\) takes the form as in \eqref{eq:JacNS}.
\end{proxexamples}

	\section{Conclusions}\label{sec:Conclusions}
		A forward-backward truncated-Newton method (\refFBTN[]) is proposed, that minimizes the sum of two convex functions one of which Lipschitz continuous and twice continuously differentiable.
Our approach is based on the forward-backward envelope (FBE), a continuously differentiable tight lower bound to the original (nonsmooth and extended-real valued) cost function sharing minima and minimizers.
The method requires forward-backward steps, Hessian evaluations of the smooth function and Clarke Jacobians of the proximal map of the nonsmooth term.
Explicit formulas of Clarke Jacobians of a wide variety of useful nonsmooth functions are collected from the literature for the reader's convenience.
The higher-order operations are needed for the computation of symmetric and positive semidefinite matrices that serve as surrogate for the Hessian of the FBE, allowing for a generalized (regularized, truncated-) Newton method for its minimization.
The algorithm exhibits global \(Q\)-linear convergence under an error bound condition, and \(Q\)-superlinear or even \(Q\)-quadratic if an additional semismoothness assumption at the limit point is satisfied.


	\clearpage
	\begin{appendix}
		\section{Auxiliary results}
			\begin{lem}\label{rem:coercive}%
	Any proper lsc convex function with nonempty and bounded set of minimizers is level bounded.
	\begin{proof}
		Let \(h\) be such function; to avoid trivialities we assume that \(\dom h\) is unbounded.
		Fix \(x_\star\in\argmin h\) and let \(R>0\) be such that \(\argmin h\subseteq B\coloneqq\ball{x_\star}{R}\).
		Since \(\dom h\) is closed, convex and unbounded, it holds that \(h\) attains a minimum on the compact set \(\boundary B\), be it \(m\), which is strictly larger than \(h(x_\star)\) (since \(\dist(\argmin h,\boundary B)>0\) due to compactness of \(\argmin h\) and openness of \(B\)).
		For \(x\notin B\), let \(s_x=x_\star+R\tfrac{x-x_\star}{\|x-x_\star\|}\) denote its projection onto \(\boundary B\), and let \(t_x\coloneqq\tfrac{\|x-x_\star\|}{R}\geq1\).
		Then,
		\begin{align*}
			h(x)
		{}={} &
			h\bigl(x_\star+t_x(s_x-x_\star))
		{}\geq{}
			h(x_\star)+t_x\bigl(h(s_x)-h(x_\star)\bigr)
		{}\geq{}
			h(x_\star)+t_x\bigl(m-h(x_\star)\bigr)
		\end{align*}
		where in the first inequality we used the fact that \(t_x\geq 1\).
		Since \(m-h(x_\star)>0\) and \(t_x\to\infty\) as \(\|x\|\to\infty\), we conclude that \(h\) is coercive, and thus level bounded.
	\end{proof}
\end{lem}

\begin{lem}\label{thm:eigs}%
	Let \(H\in\sym_+(\R^n)\) with \(\lambda_{\rm max}(H)\leq 1\).
	Then \(H-H^2\in\sym_+(\R^n)\) with
	\[
		\lambda_{\rm min}(H-H^2)
	{}={}
		\min\set{
			\lambda_{\rm min}(H)(1-\lambda_{\rm min}(H)),\,
			\lambda_{\rm max}(H)(1-\lambda_{\rm max}(H))
		}.
	\]
	\begin{proof}
		Consider the spectral decomposition \(H=\trans SDS\) for some orthogonal matrix \(S\) and diagonal \(D\).
		Then,
		\(
			H-H^2
		{}={}
			\trans S\tilde DS
		\)
		where \(\tilde D=D-D^2\).
		Apparently, \(\tilde D\) is diagonal, hence the eigenvalues of \(H-H^2\) are exactly
		\(
			\set{\lambda-\lambda^2}[\lambda\in\eigs(H)]
		\).
		The function \(\lambda\mapsto\lambda-\lambda^2\) is concave, hence the minimum in \(\eigs(\tilde H)\) is attained at one extremum, that is, either at \(\lambda=\lambda_{\rm min}(H)\) or \(\lambda=\lambda_{\rm max}(H)\), which proves the claim.
	\end{proof}
\end{lem}

\begin{lem}\label{thm:RLip}%
	For any \(\gamma\in(0,\nicefrac{2}{L_f})\) the forward-backward operator \(\T\) \eqref{eq:PG} is nonexpansive (in fact, \(\frac{2}{4-\gamma L_f}\)-averaged), and the residual \(\Res\) is Lipschitz continuous with modulus \(\frac{4}{\gamma(4-\gamma L_f)}\).
	\begin{proof}
		By combining \cite[Prop. 4.39 and Cor. 18.17]{bauschke2017convex} it follows that the \emph{gradient descent operator} \(x\mapsto \Fw x\) is \(\nicefrac{\gamma L_f}{2}\)-averaged.
		Moreover, since the proximal mapping is \(\nicefrac12\)-averaged \cite[Prop. 12.28]{bauschke2017convex} we conclude from \cite[Prop. 4.44]{bauschke2017convex} that the forward-backward operator \(T_\gamma\) is \(\alpha\)-averaged with \(\alpha=\frac{2}{4-\gamma L_f}\), thus nonexpansive \cite[Rem. 4.34(i)]{bauschke2017convex}.
		Therefore, by definition of \(\alpha\)-averagedness there exists a \(1\)-Lipschitz continuous operator \(S\) such that \(T_\gamma=(1-\alpha)\id+\alpha S\) and consequently the residual
		\(
			\Res
		{}={}
			\tfrac1\gamma
			\bigl(\id-T_\gamma\bigr)
		{}={}
			\tfrac\alpha\gamma
			(\id-S)
		\)
		is \((\nicefrac{2\alpha}{\gamma})\)-Lipschitz continuous.
	\end{proof}
\end{lem}

	\end{appendix}


	\bibliographystyle{plain}
	\bibliography{TeX/Bibliography.bib}

\end{document}